\documentclass[a4paper,10pt]{article}
\usepackage{a4, amsxtra,amsmath,amsfonts,amscd, amssymb, mathrsfs}
\usepackage[all]{xy}
\usepackage{graphicx}
\usepackage{ifpdf}
\usepackage{color}
\newtheorem{thm}{Theorem}[section]
\newtheorem{prop}[thm]{Proposition}
\newtheorem{cor}[thm]{Corollary}
\newtheorem{lem}[thm]{Lemma}
\newtheorem{ex}[thm]{Example}
\newtheorem{Def}[thm]{Definition}
\newtheorem{rem}[thm]{Remark}
\newtheorem{art}[thm]{}

\newcommand{\Hom}{{\rm Hom}}

\newcommand{\Div}{{\rm div}}
\newcommand{\cyc}{{\rm cyc}}

\newcommand{\Spec}{{\rm Spec}}

\newcommand{\id}{{\rm id}}

\newcommand{\supp}{{\rm supp}}

\newcommand{\Ccal}{{\mathscr C}}
\newcommand{\Dcal}{{\mathscr D}}

\newcommand{\Ocal}{{\mathscr O}}

\newcommand{\Ucal}{{\mathscr U}}

\newcommand{\Xcal}{{\mathscr X}}

\newcommand{\Ycal}{{\mathscr Y}}

\newcommand{\cdop}{{\mathbb C}}

\newcommand{\qdop}{{\mathbb Q}}
\newcommand{\ndop}{{\mathbb N}}
\newcommand{\rdop}{{\mathbb R}}

\newcommand{\ldop}{{\mathbb L}}

\newcommand{\adop}{{\mathbb A}}
\newcommand{\zdop}{{\mathbb Z}}

\newcommand{\proof}{\noindent {\bf Proof: \/}}
\newcommand{\qed}{{ \hfill $\square$}}

\newcommand{\Trop}{{\rm Trop}}
\newcommand{\trop}{{\rm trop}}

\newcommand{\mb}{{\mathbf m}}

\newcommand{\Sb}{{\mathbf s}}

\newcommand{\zb}{{\mathbf z}}
\newcommand{\Tor}{{\mathbb G}_m^r}

\newcommand{\Xan}{{X^{\rm an}}}
\newcommand{\Yan}{{Y^{\rm an}}}
\newcommand{\Tan}{{T^{\rm an}}}
\newcommand{\Uan}{{U^{\rm an}}}

\newcommand{\relint}{{\rm relint}}

\newcommand{\kcirc}{{ K^\circ}}
\newcommand{\ktilde}{{ \tilde{K}}}

\newcommand{\In}{{\rm in}}

\newcommand{\phitrop}{{\varphi_{\rm trop}}}

\newcommand{\del}{{\partial}}
\newcommand{\delbar}{{\bar{\partial}}}

\title{Forms and currents on the analytification of an algebraic variety (after Chambert-Loir and Ducros)}
\author{Walter Gubler}
\date{\today}
\begin{document}

\maketitle

\begin{abstract}
Chambert-Loir and Ducros have recently introduced real differential forms and currents on Berkovich spaces. In these notes, we survey this new theory and we will compare it with tropical algebraic geometry. 

\vspace{2mm}
\end{abstract}




\section{Introduction}


Antoine Chambert-Loir and Antoine Ducros have recently written the preprint ``For\-mes diff\'erentielles r\'eelles et courants sur les espaces de  Berkovich'' (see \cite{CD12}). This opens the door for applying methods from differential geometry also at non-archimedean places. We may think of possible applications for Arakelov theory or for non-archimedean dynamics. 
In the Arakelov theory developed by Gillet and Soul\'e \cite{GS90}, contributions of the $p$-adic places are described in terms of algebraic intersection theory on regular models over the valuation ring. The existence of such models  usually requires the existence of resolution of singularities which is not known in general. Another disadvantage is that canonical metrics of line bundles on abelian varieties with bad reduction can not be described in terms of models. In the case of curves, there is an analytic description of Arakelov theory also at finite places due to Chinburgh--Rumely \cite{CR93}, Thuillier \cite{Th05} and Zhang \cite{Zh93}. Now the paper of Chambert-Loir and Ducros provides us with an analytic formalism including $(p,q)$-forms, currents and differential operators $d',d''$ such that the crucial Poincar\'e--Lelong equation holds. This makes hope that we get also an analytic description of the $p$-adic contributions in Arakelov theory. In Amaury Thuillier's thesis \cite{Th05}, he has given a non-archimedean potential theory on curves. For the case of the projective line, we refer to the book of Baker and Rumely \cite{BR10} with various applications to non-archimedean dynamics. Again, we may hope to use the paper of Chambert-Loir and Ducros to give generalizations to higher dimensions.

The purpose of the present paper is to summarize the preprint \cite{CD12} and to compare it with tropical algebraic geometry. We will assume that $K$ is an algebraically closed field endowed with a (non-trivial) complete non-archimedean absolute value $|\phantom{a}|$. Let $v:=-\log|\phantom{a}|$ be the corresponding valuation and let $\Gamma:=v(K^\times)$ be the value group. Note that the residue field $\tilde{K}$ is also algebraically closed. For the sake of simplicity, we will restrict mostly to the case of an algebraic variety $X$ over $K$. In this case, there is quite an easy description of the associated analytic space $\Xan$ and so we require less knowledge about the theory of Berkovich analytic spaces than in \cite{CD12}. The main idea is quite simple: Suppose that $X$ is an $n$-dimensional closed subvariety of the split multiplicative torus $T=\Tor$. Then there is a tropicalization map $\trop:\Tan \rightarrow \rdop^r$. Roughly speaking, the map is given by applying the valuation $v$ to the coordinates of the points. Tropical geometry says that the tropical variety $\Trop(X):=\trop(\Xan)$ is a weighted polyhedral complex of pure dimension $n$ satisfying a certain balancing condition. The thesis of Lagerberg \cite{La12} gives a formalism of $(p,q)$-superforms on $\rdop^r$ together with differential operators $d',d''$ similar to $\del, \delbar$ in complex analytic geometry. Using the tropicalization map, we have a pull-back of these forms and differential operators to $\Xan$. In general, we may cover an arbitrary algebraic variety $X$ of pure dimension $n$ by very affine open charts $U$ which means that $U$ has a closed immersion to $\Tor$ and we may apply the above to define $(p,q)$-forms and currents on $\Xan$. Chambert-Loir and Ducros prove that there is an integration of compactly supported $(n,n)$-forms on $\Xan$ with the formula of Stokes and the Poincar\'e--Lelong formula. The main result of the paper \cite{CD12} is that the non-archimedean Monge-Amp\`ere measures, which were introduced by Chambert-Loir \cite{Ch06} directly as Radon measures on $\Xan$, may be written as an $n$-fold wedge product of first Chern currents. We will focus in this paper on the basics and so we will omit a description of this important result here.

\vspace{3mm}

\centerline{\it Terminology}

In $A \subset B$, $A$ may be  equal to $B$. The complement of $A$ in $B$ is denoted by $B \setminus A$ \label{setminus} as we reserve $-$ for algebraic purposes. The zero is included in $\ndop$ and in $\rdop_+$.

All occurring rings and algebras are  with $1$. If $A$ is such a ring, then the group of multiplicative units is denoted by $A^\times$. A variety over a field is an irreducible separated reduced scheme of finite type. We denote by $\overline{F}$ an algebraic closure of the field $F$. 

The terminology from convex geometry is introduced in \S \ref{Superforms and supercurrents on euclidean space} and \S  \ref{Superforms on polyhedral complexes}. Note that polytopes and polyhedra are assume to be convex. 




\vspace{3mm}
\small
This paper was written as a backup for my survey talk at the Simons Sypmosium on ``Non-Archimedean and Tropical Geometry'' in St. John from 1.4--5.4.2013. Many thanks to the organizers Matt Baker and Sam Payne for the invitation and the Simons Foundation for the support. 
The author thanks 
Antoine Chambert-Loir,  Julius Hertel, Klaus K\"unnemann, Hartwig Mayer, Jascha Smacka and Alejandro Soto for helpful comments. I am very grateful to the referee for his careful reading and his suggestions.

\normalsize

\section{Superforms and supercurrents on $\rdop^r$} \label{Superforms and supercurrents on euclidean space}

In this section, we recall the construction of superforms and supercurrents introduced by Lagerberg (see \cite{La12}, \S 2). They are real analogues of complex $(p,q)$-forms or currents on $\cdop^r$. So let us first recall briefly the definitions in complex analytic geometry. On $\cdop^r$, we have the holomorphic coordinates $z_1, \dots, z_r$. A $(p,q)$-form $\alpha$ is given by 
$$\alpha = \sum_{I,J} \alpha_{IJ} dz_I \wedge d\overline{z}_J,$$
where $I$ (resp. $J$) ranges over all subsets of $\{1,\dots,r\}$ of cardinality $p$ (resp. $q$) and where the $\alpha_{IJ}$ are smooth functions. Here, use the convenient notation $dz_I:=dz_{i_1} \wedge \cdots \wedge dz_{i_p}$ and $d\overline{z}_J:=d\overline{z}_{j_1} \wedge \cdots \wedge d\overline{z}_{j_q}$ for the elements $i_1< \dots < i_p$ of $I$ and $j_1< \dots < j_q$ of $J$. We have linear differential operators $d'$, $d''$ and $d=d'+d''$ on differential forms which are determined by the rules $$d'f=\sum_{i=1}^r\frac{\partial 
f}{\partial z_i} dz_i, \quad d''f = \sum_{j=1}^r \frac{\partial 
f}{\partial \overline{z}_j} d\overline{z}_j$$
for smooth complex functions $f$ on $\cdop^r$. Very often, these differential operators are denoted by $\del:=d'$ and $\delbar:=d''$. A current is a continuous linear functional on the space of differential forms on $\cdop^r$. Continuity is with respect to uniform convergence of finitely many derivatives on compact subsets. Differential forms may be viewed as currents using integration and the differential operators $d,d',d''$ extend to currents. For details, we refer to \cite{De12}, Chapter I, or to \cite{GH78}. 

The goal of this section is to give a real analogue in the following setting: Let $N$ be a free abelian group of rank $r$ with dual abelian group $M:=\Hom(N,\zdop)$. For convenience, we choose a basis $e_1, \dots, e_r$ of $N$ leading to coordinates $x_1, \dots , x_r$ on $N_\rdop$. Our  constructions will depend  only on the underlying real affine structure and the integration at the end will depend  on the underlying integral $\rdop$-affine structure, but not on the choice of the coordinates. Here, an {\it integral $\rdop$-affine space} is a real affine space whose underlying real vector space has an integral structure, i.e. it comes with a complete lattice. The definition of the integrals in \cite{CD12} does use calibrations which makes the integrals in some sense unnatural. In the case of an underlying canonical integral structure (which is the case for tropicalizations), there is a canonical calibration (as in \cite{CD12}, \S 3.5) and both definitions of the integrals are the same.

\begin{art} \rm \label{superforms}
Let $A^k(U,\rdop)$ be the space of smooth real differential forms of degree $k$ on an open subset $U$ of $N_\rdop$, then a {\it superform of bidegree $(p,q)$} on $U$ is an element of 
$$A^{p,q}(U):=A^p(U,\rdop) \otimes_{C^\infty(U)} A^q(U,\rdop)=C^{\infty}(U)\otimes_\zdop \Lambda^p M \otimes_\zdop \Lambda^q M.$$ Formally, such a superform $\alpha$ may be written as 
$$\alpha = \sum_{|I|=p, |J|=q}  \alpha_{IJ} d'x_I \wedge {d''x_J}$$
where $I$ (resp. $J$) consists of $i_1 < \dots < i_p$ (resp. $j_1 < \dots < j_q$), $\alpha_{IJ} \in C^\infty(U)$ and 
$$d'x_I \wedge {d''x_J} := (dx_{i_1} \wedge \dots \wedge dx_{i_p}) \otimes (dx_{j_1} \wedge \dots \wedge dx_{j_q}).$$
The wedge product is defined in the usual way on the space of superforms $A(U):= \bigoplus_{p,q \leq n} A^{p,q}(U)$ which means that $d'x_i$ and $d'x_j$ anticommute. There is a canonical $C^{\infty}(U)$-linear isomorphism $J^{p,q}:A^{p,q}(U) \rightarrow A^{q,p}(U)$ obtained by switching factors in the tensor product. The inverse of $J^{p,q}$ is $J^{q,p}$. We call $\alpha \in A^{p,p}(U)$ {\it symmetric} if $J^{p,p} \alpha = (-1)^p \alpha$. 
\end{art}

\begin{art} \rm \label{differential operators}
There is a differential operator $d':A^{p,q}(U) \rightarrow A^{p+1,q}(U)$ given by 
$$d' \alpha :=  \sum_{|I|=p, |J|=q}  \sum_{i=1}^r \frac{\partial \alpha_{IJ}}{\partial x_i} d'x_i \wedge d'x_I \wedge {d''x_J}.$$
This does not depend on the choice of coordinates as $d'=d \otimes \id$ on $A^{p,q}(U)=A^p(U, \rdop) \otimes_\zdop \Lambda^q M$ is an intrinsic characterization using the classical differential $d$ on the  space $A^p(U,\rdop)$ of real smooth $p$-forms.  Similarly, we define a differential operator 
$d'': A^{p,q}(U) \rightarrow A^{p,q+1}(U)$ by
$$d'' \alpha :=  \sum_{|I|=p, |J|=q}  \sum_{j=1}^r \frac{\partial \alpha_{IJ}}{\partial x_j} {d''x_j}\wedge d'x_I \wedge {d''x_J}.$$
By linearity, we extend these differential operators to $A(U)$. Moreover, we set $d:=d' + d''$.
\end{art}

\begin{art} \rm \label{affine pullback}
If $N'$ is  a free abelian group of rank $r'$ and if $F:N'_\rdop \rightarrow N_\rdop$ is an affine map with $F(V) \subset U$ for an open subset $V$ of $N'_\rdop$, then we have a well-defined pull-back $F^*:A^{p,q}(U) \rightarrow A^{p,q}(V)$ given as usual. The affine pull-back commutes with the differential operators $d$, $d'$ and $d''$. The pull-back is defined more generally for smooth maps, but then it does not necessarily commute with $d$, $d'$ and $d''$.
\end{art}

\begin{art} \rm \label{integration}
Let $A_c(U)$ denote the space of superforms on $U$ with compact support in $U$. Recall that $r$ is the rank of $M$. For $\alpha \in  A_c(U)$, we define 
$$ \int_U \alpha := (-1)^{\frac{r(r-1)}{2}}\int_U \alpha_{LL} dx_1 \wedge \dots \wedge dx_r$$
with $L=\{1,\dots, r\}$ and the usual integration of $r$-forms with respect to the orientation induced by the choice of coordinates on the right hand side. If $F$ is an affine map as in \ref{affine pullback} and if $r=r'$, then we have the transformation formula
\begin{equation} \label{trafo1}
\int_V F^*(\alpha) = |\det(F)| \int_V \alpha
\end{equation}
(see \cite{La12}, equation (2.3)). We conclude that the definition of the integral depends only on the underlying integral $\rdop$-affine structure of $N_\rdop$.

The sign $ (-1)^{\frac{r(r-1)}{2}}$ is explained by the fact that we want $d'x_1 \wedge d''x_1 \wedge \dots \wedge d'x_r \wedge d''x_r$ to be a positive $(r,r)$-superform and hence
$$\int_U f d'x_1 \wedge d''x_1 \wedge \dots \wedge d'x_r \wedge d''x_r =
\int_U f dx_1 \wedge \dots \wedge dx_r$$ 
for any $f \in C^\infty_c(U)$ (see \cite{CD12} for more details about positive forms).
\end{art}

\begin{art} \rm \label{polyhedra}
Now let $\sigma$ be a {\it polyhedron} of dimension $n$ in $N_\rdop$. By definition, $\sigma$ is the intersection of finitely many halfspaces $H_i:=\{\omega \in N_\rdop \mid \langle u_i , \omega \rangle \leq c_i\}$ with $u_i \in M_\rdop$ and $c_i \in \rdop$.  A {\it polytope} is a bounded polyhedron. We say that $\sigma$ is an {\it integral $G$-affine polyhedron} for a subgroup $G$ of $\rdop$ if we may choose all $u_i \in M$ and all $c_i \in G$. In this case,  we have a canonical integral $\rdop$-affine structure on the affine space $\adop_\sigma$ generated by $\sigma$. If $\ldop_\sigma$ is the underlying real vector space of $\adop_\sigma$, then this integral structure is given by the lattice $N_\sigma:=\ldop_\sigma \cap N$. Using \ref{affine pullback} and the above, we get a well-defined integral $\int_\sigma \alpha$ for any $\alpha \in A_c^{n,n}(U)$, where $U$ is an open neighbourhood of $\sigma$. 
\end{art}

\begin{art} \label{contraction} \rm 
In \cite{CD12}, integration is described in terms of a contraction: Similarly as in differential geometry, we may view a superform $\alpha \in A^{p,q}(U)$ as a multilinear map $$N_\rdop^{p+q} \longrightarrow C^{\infty}(U), \quad (n_1,\dots, n_{p+q}) \mapsto \alpha(n_1, \dots, n_{p+q})$$ which is alternating in the variables $(n_1, \dots,n_p)$ and also in  $(n_{p+1}, \dots, n_{p+q})$. Let $I \subset \{1,\dots, p+q\}$ be a subset of cardinality $s$ with $s'$ elements contained in $\{1,\dots, p\}$ and hence $s''=s-s'$ elements in $\{p+1, \dots,p+q\}$.  Given  vectors $v_1, \dots, v_s \in N_\rdop$, the {\it contraction} $\langle \alpha;v_1, \dots, v_s \rangle_I \in A^{p-s',q-s''}(U)$ is given by inserting $v_1, \dots, v_s$ for the variables $(n_i)_{i \in I}$ of the above multilinear function. 

Using the basis $e_1, \dots, e_r$ of $N$ and assuming $\alpha \in A_c^{r,r}(U)$, the contraction $ \langle \alpha; e_1, \dots, e_r \rangle_{\{r+1,\dots, 2r\}}$ is a $(r,0)$-superform which may be viewed as a classical $r$-form on $U$. Then it is immediately clear from the definitions that we have
$$\int_U \alpha = (-1)^{\frac{r(r-1)}{2}} \int_U \langle \alpha; e_1, \dots, e_r \rangle_{\{r+1,\dots, 2r\}}$$ 
where we use the usual integration of $r$-forms on the right. Of course, there is no preference to contract with respect to the last $r$ variables. Similarly, may view $\langle \alpha; e_1, \dots, e_r \rangle_{\{1,\dots, r\}} \in A_c^{0,r}(U)$ as a classical $r$-form and we have
$$\int_U \alpha = (-1)^{\frac{r(r-1)}{2}} \int_U \langle \alpha; e_1, \dots, e_r \rangle_{\{1,\dots, r\}}.$$
\end{art}

Next, we are looking for an analogue of Stokes' theorem for superforms. 

\begin{art} \rm \label{hyperplane integration} 
Let $H$ be an integral $\rdop$-affine halfspace in $N_\rdop$. This means that $H=\{\omega \in N_\rdop \mid \langle u, \omega \rangle \leq c\}$ for some $u \in M$ and  $c \in \rdop$.   Using a translation, we may assume that $c=0$ and hence the boundary $\partial H$ is a linear subspace of $N_\rdop$. Let $[\omega_{\partial H, H}]$ be the generator of $N/(N \cap \partial H)\cong \zdop$ which points outwards, i.e. there is $u_{\partial H, H} \in M$ such that $u_{\partial H, H}(H) \leq 0$ and $u_{\partial H, H}(\omega_{\partial H, H})=1$. We choose a representative $\omega_{\partial H, H} \in N$ and we note also that $u_{\partial H, H}$ is uniquely determined by the above properties. 
\end{art}

\begin{art} \rm \label{setup for Stokes for polyhedrons}
Let $U$ be an open subset of $N_\rdop$ and let $\sigma$ be an $r$-dimensional integral $\rdop$-affine polyhedron contained in $U$. For any closed face $\rho$ of codimension $1$, let $\omega_{\rho,\sigma}:=\omega_{\partial H, H}$ using \ref{hyperplane integration}   for the affine hyperplane $\partial H$ generated by $\rho$ and the corresponding halfspace containing $\sigma$. We note that $\omega_{\rho,\sigma} \in N$ is determined up to addition with elements in $N_\rho=N \cap \ldop_\rho$, where $\ldop_\rho$ is the linear hyperplane parallel to $\rho$. 

For $\eta \in A_c^{r-1,r}(U)$, we have introduced the contraction 
$\langle \eta; \omega_{\rho,\sigma} \rangle_{\{2r-1\}}$ as an element of  $A_c^{r-1,r-1}(U)$ which is obtained by inserting the vector $\omega_{\rho,\sigma}$ for the $(2r-1)$-th argument of the corresponding multilinear function (see \ref{contraction}). Note that the restriction of this contraction to $\rho$ does not depend on the choice of the representative $\omega_{\rho,\sigma}$. Then we define
$$\int_{\partial \sigma} \eta := \sum_\rho \int_\rho \langle \eta; \omega_{\rho,\sigma} \rangle_{\{2r-1\}},$$
where $\rho$ ranges over all closed faces of $\sigma$ of codimension $1$. On the right, we use the integrals of $(r-1,r-1)$-superforms from \ref{integration}. For $\eta \in A_c^{r,r-1}(U)$, we define similarly \footnote{There is a sign mistake in the published version which is here corrected in \textcolor{red}{red}.}
$$\int_{\partial \sigma} \eta := \textcolor{red}{-}\sum_\rho \int_\rho \langle \eta; \omega_{\rho,\sigma} \rangle_{\{\textcolor{red}{1}\}} .$$
Note that the integrals do depend only on the integral $\rdop$-affine structure of $N_\rdop$ but do not depend on the choice of the orientation of $N_\rdop$.

If $\sigma$ is an integral $\rdop$-affine polyhedron of any dimension $n$ and if  $\eta \in A_c^{n-1,n}(U)$ for an open subset $U$ of $N_\rdop$ containing $\sigma$, then we define $\int_{\partial \sigma} \eta$ by applying the above to the affine space $\adop_\sigma$ generated by $\sigma$ and to the pull-back of $\eta$ to $\adop_\sigma$. We give now a concrete description  of $\int_{\partial \sigma} \eta$ in terms of integrals over classical $(n-1)$-forms. For every closed face $\rho$ of $\sigma$, let $N_\sigma = \ldop_\sigma \cap N$ be the canonical integral structure on the affine space generated by $\sigma$. If  $e_1^\rho, \dots, e_{n-1}^\rho$ is a basis of $N_\rho$, then $\omega_{\rho,\sigma}, e_1^\rho, \dots, e_{n-1}^\rho $ is a basis of $N_\sigma$. We note that the contraction  $ \langle \eta; \omega_{\rho,\sigma}, e_1^\rho, \dots, e_{n-1}^\rho \rangle_{\{n, \dots, 2n-1\}}$ may be viewed as a classical $(n-1)$-form on $U$ and hence we get
$$\int_{\partial \sigma} \eta = \sum_\rho \int_\rho \langle \eta; \omega_{\rho,\sigma} \rangle_{\{2n-1\}}=  (-1)^{\frac{n(n-1)}{2}} \sum_\rho \int_\rho \langle \eta; \omega_{\rho,\sigma}, e_1^\rho, \dots, e_{n-1}^\rho \rangle_{\{n, \dots, 2n-1\}}.$$
\end{art}


\begin{prop}[Stokes' formula] \label{Stokes for polyhedrons}
Let $\sigma$ be an $n$-dimensional integral $\rdop$-affine polyhedron contained in the open subset $U$ of $N_\rdop$. For any $\eta' \in A_c^{n-1,n}(U)$ and any $\eta'' \in A_c^{n,n-1}(U)$, we have
$$\int_\sigma d' \eta' = \int_{\partial \sigma}  \eta', \quad
\int_\sigma  d'' \eta'' =  \int_{\partial \sigma} \eta''.$$
\end{prop}

\proof This is just a reformulation of \cite{La12}, Proposition 2.3, in the case of a polyhedron using the formalism introduced above. In the quoted result, the boundary was assumed to be smooth, but as the classical Stokes' formula holds also for polyhedra (see \cite{Wa83}, 4.7), this applies here as well. \qed

\begin{prop}[Green's formula] \label{Green for polyhedrons}
We consider  an $n$-dimensional integral $\rdop$-affine polyhedron $\sigma$ contained in the open subset $U$ of $N_\rdop$. Assume that $\alpha \in A^{p,p}(U)$ and $\beta \in A^{q,q}(U)$ are symmetric with $p+q=n-1$ and that the intersection of the supports of $\alpha$ and $\beta$ is compact. Then we have
$$\int_\sigma \alpha \wedge d'd''\beta - \beta \wedge d'd''\alpha = \int_{\partial \sigma} \alpha \wedge d''\beta  - \beta \wedge d'' \alpha.$$
\end{prop}

\proof This follows from Stokes' formula as in \cite{CD12}, Lemma 1.3.8. \qed

\begin{art} \rm \label{supercurrents}
A {\it supercurrent} on $U$ is a continuous linear functional on $A_c^{p,q}(U)$ where the latter is a locally convex vector space in a similar way as in the classical case. We denote the space of such supercurrents by $D_{p,q}(U)$. As usual , we  define the linear differential operators $d$, $d'$ and $d''$ on $D(U):=\bigoplus_{p,q}D_{p,q}(U)$ by using $(-1)^{p+q+1}$ times the dual of the corresponding differential operator on $A_c^{p,q}(U)$. The sign is chosen in such a way that the canonical embedding $A^{p,q}(U) \rightarrow D_{r-p,r-q}(U)$ is compatible with the operators $d$, $d'$ and $d''$. Here, $\alpha \in A^{p,q}(U)$ is mapped to $[\alpha] \in D_{r-p,r-q}(U)$ given by $[\alpha](\beta)=\int_{N_\rdop} \alpha \wedge \beta$ for any $\beta \in A_c^{r-p, r-q}(U)$. 
\end{art}

\section{Superforms on polyhedral complexes} \label{Superforms on polyhedral complexes}

We keep the notions from the previous section and we will extend them to the setting of polyhedral complexes. We will introduce tropical cycles and we will characterize them as closed currents of integrations over weighted integral $\rdop$-affine polyhedral complexes. 

\begin{art} \rm \label{polyhedral complexes}
A {\it polyhedral complex} $\Ccal$ in $N_\rdop$ is a finite set of polyhedra with the following two properties: Every polyhedron in $\Ccal$ has all its closed faces in $\Ccal$. If $\Delta, \sigma \in {\Ccal}$, then $\Delta \cap \sigma$ is a closed face of $\Delta$ and $\sigma$. Note here that the empty set and also $\sigma$ are allowed as closed faces of a polyhedron $\sigma$ (see \cite{Gu12}, Appendix A, for details).

A  polyhedral complex $\Ccal$ is called {\it integral $G$-affine} for a subgroup $G$ of $\rdop$ if every polyhedron of $\Ccal$ is integral $G$-affine. The {\it support} $|\Ccal|$ of $\Ccal$ is the union of all polyhedra in $\Ccal$. The polyhedral complex $\Ccal$ is called {\it pure dimensional of dimension $n$} if every maximal polyhedron in $\Ccal$ has dimension $n$. We will often use the notation $\Ccal_k:=\{\sigma \in \Ccal \mid \dim(\sigma)=k\}$ for $k \in \ndop$.
\end{art}




\begin{art} \rm \label{tropical superforms}
Let $\Ccal$ be a polyhedral complex in $N_\rdop$. A {\it superform} on $\Ccal$ is the restriction of a superform on (an open subset of) $N_\rdop$ to $|\Ccal|$. This means that two superforms agree if their restrictions to any polyhedron of $|\Ccal|$ agree. Let $A(\Ccal)$ be the space of superforms on $\Ccal$. It is an alternating algebra with respect to the induced wedge product. We have also differential operators $d$, $d'$ and $d''$ on $A(\Ccal)$ given by restriction of the corresponding operators on $A(N_\rdop)$. Let $A^{p,q}(\Ccal)$ be the space of $(p,q)$-superforms on $\Ccal$. The {\it support} of $\alpha \in A(\Ccal)$ is the complement of $\{\omega \in |\Ccal| \mid \text{$\alpha$ vanishes identically in a neighbourhood of $\omega$}\}$ in $|\Ccal|$. We denote by $A_c^{p,q}(\Ccal)$ the subspace of $A^{p,q}(\Ccal)$ of superforms of compact support.

Let  $N'$ be  a free abelian group of rank $r'$ and let $F:N'_\rdop \rightarrow N_\rdop$ be an affine map. Suppose that $\Ccal'$ is a polyhedral complex of $N'_\rdop$ with $F(|\Ccal'|) \subset |\Ccal|$, then the pull-back  in \ref{affine pullback} induces a pull-back $F^*:A^{p,q}(\Ccal) \rightarrow A^{p,q}(\Ccal')$.
\end{art}

\begin{art} \rm \label{weights}
A polyhedral complex $\Dcal$ {\it subdivides} the polyhedral complex $\Ccal$ if they have the same support and if every polyhedron $\Delta$ of $\Dcal$ is contained in a polyhedron of $\Ccal$. In this case, we say that $\Dcal$ is a {\it subdivision} of $\Ccal$. All our constructions here will be compatible with subdivisions. This is no problem for the definition of superforms on $\Ccal$ as they depend only on the support $|\Ccal|$.

A {\it weight} on a pure dimensional polyhedral complex $\Ccal$ is a function $m$ which assigns to every maximal polyhedron $\sigma \in \Ccal$ a number $m_\sigma \in \zdop$. Then we get a canonical weight on every subdivision of $\Ccal$. For a weighted polyhedral complex $(\Ccal, m)$, only the polyhedra $\Delta \in \Ccal$ which are contained in a maximal dimensional $\sigma \in \Ccal$ with $m_\sigma \neq 0$ are of interest. They form a subcomplex $\Dcal$ of $\Ccal$ and we define the support of $(\Ccal,m)$ as the support of $\Dcal$. The polyhedra of $\Ccal \setminus \Dcal$ will usually be neglected.
\end{art}

\begin{art} \label{tropical integration} \rm  
Let $(\Ccal,m)$ be a weighted integral $\rdop$-affine polyhedral complex of pure dimension $n$. For $\alpha \in A_c^{n,n}(\Ccal)$,  we set 
$$\int_{(\Ccal,m)} \alpha := \sum_{\sigma \in \Ccal_n} m_\sigma \int_\sigma \alpha,$$
where we use integration from \ref{integration} on the right. We define integrals over the boundary of $\Ccal$ for a superform $\beta$ in  $A_c^{n-1,n}(\Ccal)$ or in  $A_c^{n,n-1}(\Ccal)$ by
$$\int_{\partial(\Ccal,m)} \beta = \sum_{\sigma \in \Ccal_n} m_\sigma \int_{\partial \sigma} \beta,$$
where we use the boundary integrals from \ref{setup for Stokes for polyhedrons} on the right. Note that the boundary $\partial \Ccal$  may be defined as the subcomplex consisting of the polyhedra of dimension at most $n-1$, but there is no canonical weight on $\partial \Ccal$. Indeed, the boundary integral $\int_{\partial(\Ccal,m)} \beta$ depends on the relative situation $\partial \Ccal \subset \Ccal$ because of the weight $m_\sigma$ and the contraction with respect to the vectors $\omega_{\rho,\sigma}$ used in the definitions. This is similar to the situation in real analysis where boundary integrals depend on the relative orientation. These classical boundary integrals do depend only on the restriction of the differential form to the boundary which is clearly wrong for our boundary integrals. However, it is still true that $\int_{\partial(\Ccal,m)} \beta =0$ if the support of $\beta$ is disjoint from $\partial \Ccal$.
\end{art}

\begin{prop}[Stokes' formula] \label{tropical Stokes}
Let $(\Ccal,m)$ be a weighted integral $\rdop$-affine polyhedral complex of pure dimension $n$. For any $\eta' \in A_c^{n-1,n}(\Ccal)$ and any $\eta'' \in A_c^{n,n-1}(\Ccal)$, we have
$$\int_{(\Ccal,m)} d' \eta' = \int_{\partial (\Ccal,m)}  \eta', \quad
\int_{(\Ccal,m)}  d'' \eta'' =  \int_{\partial(\Ccal,m)} \eta'.$$
\end{prop}

\proof This follows immediately from Stokes' formula for polyhedra given in Proposition \ref{Stokes for polyhedrons}.
\qed


\begin{ex} \rm \label{Dirac current}
If $(\Ccal,m)$ is a weighted integral $\rdop$-affine polyhedral complex of pure dimension $n$, then we get a  supercurrent $\delta_{(\Ccal,m)} \in D_{n,n}(N_\rdop)$ by setting $\delta_{(\Ccal,m)}(\eta)=\int_{(\Ccal,m)} \eta$ for any $\eta \in A_c^{n,n}(N_\rdop)$. 
\end{ex}

\begin{art} \rm \label{tropical cycle}
A weighted integral $\rdop$-affine polyhedral complex $(\Ccal,m)$ of pure dimension $n$  is called a {\it tropical cycle} if its weight $m$ satisfies the following {\it balancing condition}: For every $(n-1)$-dimensional $\rho \in \Ccal$, we have 
$$\sum_{\sigma \in \Ccal_n, \,\sigma \supset \rho} m_\sigma \omega_{\rho, \sigma} \in N_\rho.$$
Here, $N_\rho$ is the canonical lattice contained in the affine space generated by $\rho$ and $\omega_{\rho, \sigma} \in N_\sigma$ is the lattice vector pointing outwards of $\sigma$ (see  \ref{setup for Stokes for polyhedrons}). Tropical cycles are the basic objects in tropical geometry.
\end{art} 

\begin{prop} \label{tropical cycle iff closed}
Let $(\Ccal,m)$ be a weighted  integral $\rdop$-affine polyhedral complex of pure dimension $n$  on $N_\rdop$. Then the following conditions are equivalent:
\begin{itemize}
\item[(a)] $(\Ccal,m)$ is a tropical cycle;
\item[(b)] $\delta_{(\Ccal,m)}$ is a $d'$-closed supercurrent on $N_\rdop$;
\item[(c)] $\delta_{(\Ccal,m)}$ is a $d''$-closed supercurrent on $N_\rdop$.
\end{itemize}
\end{prop}

\proof Let $\alpha \in A_c^{n-1,n}(N_\rdop)$. By Stokes' formula in Proposition \ref{tropical Stokes}, we have
\begin{equation*} \label{integration identity}
\delta_{(\Ccal,m)} (d'\alpha) = \int_{\partial (\Ccal,m)} \alpha =
\sum_\rho  \int_\rho \langle \alpha; \sum_{\sigma \supset \rho} m_\sigma \omega_{\rho,\sigma} \rangle_{\{2n-1\}},
\end{equation*}
where $\rho$ (resp. $\sigma$) ranges over all elements of $\Ccal$ of dimension $n-1$ (resp. $n$). Suppose now that $\sum_{\sigma \supset \rho} m_\sigma \omega_{\rho,\sigma} \in N_\rho$ for some $(n-1)$-dimensional $\rho \in \Ccal$. Recall that we may view $\alpha$ as a multilinear map  $N_{\rdop}^{2n-1} \rightarrow C^\infty(N_\rdop)$ which is alternating in the first $n-1$ arguments and also alternating in the last $n$ arguments. But an alternating $n$-linear map on a vector space of dimension $n-1$ is zero and hence the restriction of  $\langle \alpha; \sum_{\sigma \supset \rho} m_\sigma \omega_{\rho,\sigma} \rangle_{\{2n-1\}}$ to $\rho$ is zero. Then the above display proves (a) $\Rightarrow$ (b).

Conversely, if $\sum_{\sigma \supset \rho} m_\sigma \omega_{\rho,\sigma} \not \in N_\rho$ for some $(n-1)$-dimensional  $\rho \in \Ccal$ , then there is an $\alpha \in A_c^{n-1,n}(N_\rdop)$ such that the restriction of  $\langle \alpha; \sum_{\sigma \supset \rho} m_\sigma \omega_{\rho,\sigma} \rangle_{\{2n-1\}}$ to $\rho$ is non-zero. We may also assume that the support of $\alpha$ is disjoint from all other $(n-1)$-dimensional polyhedra of $\Ccal$. Then the above display proves (b) $\Rightarrow$ (a). The equivalence of (a) and (c) is shown similarly.
\qed

\begin{art} \label{push-forward of tropical cycles} \rm
Now let $F:N_\rdop' \rightarrow N_\rdop$ be an affine map whose underlying linear map is integral, i.e. induced by a homomorphism $N' \rightarrow N$. We will define the push-forward of a weighted  integral $\rdop$-affine polyhedral complex $(\Ccal',m)$ of pure dimension $n$  on $N_\rdop'$. 
 For details, we refer to \cite{AR10}, \S 7.  After a subdivision of $\Ccal'$, we may assume that 
$$F_*(\Ccal'):=\{F(\sigma') \mid \text{$\sigma'$ is a face of $\nu' \in \Ccal'$ with $\dim(F(\nu'))=n$}\}$$
is a polyhedral complex in $N_\rdop$. We define the multiplicity of an $n$-dimensional $F(\sigma') \in F_*(\Ccal')$ by
$$m_{F(\sigma')}:=\sum_{\nu' \in \Ccal'_n,\,\nu' \subset F^{-1}(F(\sigma'))}[M_{\nu'}':M_{F(\sigma')}] m_{\nu'} .$$
Endowed with these multiplicities, we get a weighted integral $\rdop$-affine polyhedral complex $F_*(\Ccal',m)$ of $N_\rdop$. If $(\Ccal',m)$ is a tropical cycle, then $F_*(\Ccal',m)$ is also a tropical cycle.  It might happen that $F_*(\Ccal',m)$ is empty, then we get the tropical zero cycle.  
\end{art}


\begin{prop}[projection formula] \label{projection formula for tropical superforms}
Using the assumptions above and $\alpha \in A_c^{n,n}(F_*(\Ccal'))$, we have $\int_{F_*(\Ccal',m)} \alpha = \int_{(\Ccal',m)} F^*(\alpha)$.
\end{prop}

\proof Let $\sigma'$ be an $n$-dimensional polyhedron of $\Ccal'$. Then $\sigma := F(\sigma')$ is an integral $\rdop$-affine polyhedron in $N_\rdop$. We assume for the moment that $\sigma$ is also $n$-dimensional. As above, we consider the lattice $N_\sigma := N \cap \ldop_\sigma$ in $N_\rdop$, where $\ldop_\sigma$ is the linear space which is a translate of the affine space generated by $\sigma$. Let $A$ be the matrix of the homomorphism $F:N_{\sigma'} \rightarrow N_\sigma$ with respect to integral bases. Then we have $|\det(A)|=[N_{\sigma'}:N_\sigma]$ and hence the transformation formula \eqref{trafo1} shows 
\begin{equation} \label{trafo2}
 \int_{\sigma'} F^*\alpha =[N_{\sigma'}:N_\sigma]  \int_\sigma \alpha.
\end{equation}
If $\dim(\sigma)<n$, then both sides are zero and hence  formula \eqref{trafo2} is true in any case. Using the weighted sum over all $\sigma'$, the claim follows immediately from \eqref{trafo2}. \qed


\section{Moment maps and tropical charts}

A complex manifold is locally defined using analytic charts $\varphi:U \rightarrow \cdop^r$. The charts help to transport the analysis from $\cdop^r$ to $M$.  The idea in the non-archimedean setting is similar replacing the above charts by algebraic moment maps $\varphi:U \rightarrow \Tor$ to multiplicative tori  and the corresponding tropicalizations $\phitrop: U \rightarrow \rdop^r$. The restriction of $\varphi^{\rm an}$ to the preimage of an open analytic subset will be called a tropical chart. 

In this section, $K$ is an algebraically closed field endowed with a complete non-trivial non-archimedean absolute value $|\phantom{a}|$. Note that the residue field $\ktilde$ is also algebraically closed. Let $v:= - \log |\phantom{a}|$ be the associated valuation and let $\Gamma:=v(K^\times)$ be the value group. We will study analytifications, tropicalizations and moment maps of the  algebraic variety $X$ over $K$. This will be used in the next section to define $(p,q)$-forms on $\Xan$.

\begin{art} \rm \label{analytification}
We recall first the construction of the analytification of $X$. Let $U = \Spec(A)$ be an open affine subset of $X$, then $\Uan$ is the set of multiplicative seminorms on $A$ extending the given absolute value $|\phantom{a}|$ on $K$. This set is endowed with the topology generated by the functions $\Uan \rightarrow \rdop, p \mapsto p(a)$ with $a$ ranging over $A$. By glueing, we get a topological space $\Xan$ which is connected locally compact and Hausdorff. 
We can endow it with a sheaf of analytic functions leading to a Berkovich analytic space over $K$ which we call the {\it analytification} of $X$. For a morphism $\varphi:Y \rightarrow X$ of algebraic varieties over $K$, we get an analytic morphism $\varphi^{\rm an}: \Yan \rightarrow \Xan$ induced by composing the multiplicative semiorms with $\varphi^\sharp$ on suitable affine open subsets. We refer to \cite{Be90} for details,  or to \cite{BPS11}, \S 1.2, for a neat description of the analytification.
\end{art}

\begin{art} \rm \label{local dimension}
We will define some local invariants in $x \in \Xan$. On an open affine neighbourhood $U = \Spec(A)$, the point $x$ is given by a multiplicative seminorm $p$ on $A$ and we often write $|f(x)|:=p(f)$ for $f \in A$. 
Dividing out the prime ideal $I:=\{f \in A \mid p(f)=0\}$, we get a multiplicative norm on the integral domain $ B:=A/I$ which extends to an absolute value $|\phantom{a}|_x$ on the quotient field of $B$. The completion of this field is denoted by ${\mathscr H}(x)$. It does not depend on the choice of $U$ and it may be also constructed analytically. The absolute value of ${\mathscr H}(x)$ is denoted by $|\phantom{a}|$ as it extends the given absolute value on $K$. Note that the completed residue field  ${\mathscr H}(x)$ of $x$ remains the same if we replace the ambient variety $X$ by the Zariski closure of $x$ in $X$.

Let $s(x)$ be the transcendence degree of the residue field of ${\mathscr H}(x)$ over $\ktilde$. The quotient of the value group of ${\mathscr H}(x)$ by $\Gamma$ is a finitely generated abelian group and we denote its $\qdop$-rank by $t(x)$. Finally, we set $d(x):=s(x)+t(x)$. Note that Abhyankar's inequality shows that $d(x)$ is bounded by the transcendence degree of ${\mathscr H}(x)/K$. By \cite{Be90}, Proposition 9.1.3, we have $\dim(X)=\dim(V)=\sup_{x \in V}d(x)$ for every open subset $V$ of $\Xan$. 
\end{art}

\begin{ex} \rm \label{torus}
Let $T=\Tor$ be the split multiplicative torus of rank $r$ with coordinates $z_1, \dots, z_r$. Then a point $x$ of $\Tan$ could be visualized by the coordinates $z_1(x), \dots, z_r(x) \in {\mathscr H}(x)$ and the multiplicative seminorm corresponding to $x$ is given by $|f(x)|=|f(z_1(x), \dots, z_r(x))|$ for every Laurent polynomial $f$ on $T$. Conversely, every field extension $L/K$ with an absolute value extending the given absolute value on $K$ and every $(\beta_1, \dots, \beta_r) \in (L^\times)^r$ give rise to a point $x \in \Tan$ by $|f(x)|:=|f(\beta_1, \dots, \beta_r)|$. Note that $L$ and $(\beta_1, \dots, \beta_r)$ are not uniquely determined by $x$. 

In particular, we get an inclusion of $T(K)$ into $\Tan$. For every $x \in T(K)$, we have $d(x)=0$. However, there can be also other points with $d(x)=0$. If $T = {\mathbb G}_m^1$, then precisely the points of type 1 (i.e. the $K$-rational points) and the points of type 4 satisfy $d(x)=0$ (see \cite{Be90}, 1.4.4). 

Returning to the case $T = \Tor$, there are some distinguished points of $\Tan$ which behave completely different than $K$-rational points. For positive real numbers $s_1, \dots, s_r$, we define the {\it associated weighted Gauss norm} on $K[T]$ by
$$|f|_\Sb := \max_{\mb \in \zdop^r} |\alpha_\mb| \Sb^\mb$$
for every Laurent polynomial $f=\sum_{\mb \in \zdop^r} \alpha_\mb \zb^\mb \in K[T]=K[z_1^{\pm 1}, \dots, z_r^{\pm 1}]$. It follows from the Gauss Lemma that the weighted Gauss norm is a multiplicative seminorm giving rise to a point $\eta_\Sb \in \Tan$. 
The set $S(\Tan):=\{\eta_\Sb \mid s_1 >0 , \dots ,s_r >0\}$ is called the {\it skeleton} of $\Tan$. Every point $\eta_\Sb \in S(\Tan)$ satisfies $d(\eta_\Sb)=r$ (see \cite{Du12}, (0.12) and (0.13)). 
\end{ex}

\begin{art} \rm \label{tropicalization map}
Let $T:=\Tor$ be a split  multiplicative torus over $K$ with coordinates $z_1, \dots, z_r$. Then we  have the {\it tropicalization map} 
$${\trop}: \Tan \rightarrow \rdop^r, \quad p \mapsto (-\log p(z_1), \dots , -\log p(z_r)).$$
It is immediate from the definitions that the map ${\trop}$ is continuous and proper. To get a coordinate free approach, we could  use the character group $M$ and its dual $N$. Then $\trop$ is a map from $\Tan$ to $N_\rdop$. We refer to \cite{Gu12} for details about tropical geometry.
\end{art}

\begin{rem} \rm \label{tropicalization as a retraction}
Note that we have a natural section $\rdop^r \rightarrow \Tan$ of the tropicalization map. It is given by mapping the point $\omega \in \rdop^r$ to the weighted Gauss norm $\eta_\Sb$ associated to $\Sb:=(e^{-\omega_1}, \dots, e^{-\omega_r})$. It follows from \cite{Be90}, Example 5.2.12, that this section is a homeomorphism of $\rdop^r$ onto a closed subset of $\Tan$ which is the skeleton $S(\Tan)$ introduced in \ref{torus}. In this way, we may view the tropicalization map as a map from $\Tan$ onto  $S(\Xan)$.  Then it is shown in \cite{Be90}, \S 6.3, that the tropicalization map is a strong deformation  retraction of $\Tan$ onto the skeleton $S(\Tan)$. This point of view is used very rarely in our paper.
\end{rem}

\begin{art} \rm \label{tropical variety}
For a closed subvariety $Y$ of $T$ of dimension $n$, the {\it tropical variety associated to $X$} is defined by $\Trop(Y):={\trop}(\Yan)$. The Bieri--Groves theorem says that $\Trop(Y)$ is a finite union of $n$-dimensional integral $\Gamma$-affine polyhedra in $\rdop^r$. It is shown in tropical geometry that $\Trop(Y)$ is an  integral $\Gamma$-affine polyhedral complex. The polyhedral structure is only determined up to subdivision which does not matter for our constructions. We will see below that the tropical variety is endowed with a positive canonical weight $m$ satisfying the balancing condition from \ref{tropical cycle}. We get a tropical cycle of pure dimension $n$ which we  also denote by $\Trop(Y)$ forgetting the weight $m$ in the notation. 
\end{art}

\begin{art} \rm \label{tropical weights}
The {\it tropical weight} $m$ on an $n$-dimensional polyhedron $\sigma$ of $\Trop(Y)$ is defined in the following way. By density of the value group $\Gamma$ in $\rdop$, there is $\omega \in \Gamma^r \cap \relint(\sigma)$. We choose $t \in \Tor(K)$ with $\trop(t)=\omega$. Then the closure of $t^{-1}Y$ in $({\Tor})_{\kcirc}$ is a flat variety over $\kcirc$ whose special fibre is called the {\it initial degeneration} $\In_{\omega}(Y)$ of $Y$ at $\omega$. Note that $\In_{\omega}(Y)$ is a closed subscheme of $(\Tor)_{\ktilde}$. Let $m_W$ be the multiplicity of an irreducible component $W$ of $\In_{\omega}(Y)$. Then the {\it tropical weight} $m_\sigma$ is defined by $m_\sigma:=\sum_W m_W$, where $W$ ranges over all irreducible components of $\In_{\omega}(Y)$. One can show that the definition is independent of the choices of $\omega$ and $t$. It is a non-trivial fact from tropical geometry that $(\Trop(Y),m)$ is a tropical cycle (see \cite{Gu12}, \S 13, for details).
\end{art}

\begin{art} \rm \label{moment map}
For an open subset $U$ of the algebraic variety $X$, a  {\it moment map} is a morphism $\varphi:U \rightarrow T$ to a split multiplicative torus $T:=\Tor$ over $K$. The {\it tropicalization} of $\varphi$ is 
$$\phitrop:=\trop \circ \varphi^{\rm an}: \Uan \overset{\varphi^{\rm an}}{\longrightarrow} \Tan \overset{\trop}{\longrightarrow} \rdop^r.$$
Obviously, this is a continuous map with respect to the topology on the analytification $\Uan$. Note that our moment maps are algebraic  which differs from the moment maps  in \cite{CD12} which are defined analytically. 

We say that the moment map $\varphi':U' \rightarrow T'$ of the open subset $U'$ of $X$ {\it refines} the moment map $\varphi:U \rightarrow T$ if $U' \subset U$ and if there is an affine homomorphism $\psi:T' \rightarrow T$ of the multiplicative tori such that $\varphi = \psi \circ \varphi'$ on $U'$. Here, an {\it affine homomorphism} means a group homomorphism composed with a (multiplicative) translation on $T$. This group homomorphism induces a homomorphism $M \rightarrow M'$ of character lattices. Its dual is the linear part of an integral affine map $\Trop(\psi):N'_\rdop \rightarrow N_\rdop$ such that $\phitrop = \Trop(\psi) \circ \varphi_{\rm trop}'$ on $(U')^{\rm an}$.

If $\varphi_i: U_i \rightarrow T_i$ are finitely many moment maps of non-empty open subsets $U_i$ of $X$ with $i \in I$, then $U:= \bigcap_i U_i$ is a non-empty open subset of $X$ and $\varphi:U \rightarrow \prod_i T_i, x \mapsto (\varphi_i(x))_{i \in I}$ is a moment map which refines every $\varphi_i$. Moreover, it follows easily from the universal property of the product that every moment map $\varphi':U' \rightarrow T'$ which refines every $\varphi_i$ refines also $\varphi$.
\end{art}

\begin{lem} \label{density and tropicalization}
Let $\varphi:U \rightarrow \Tor$ be a moment map on an open subset $U$ of $X$ and let $U'$ be a non-empty open subset of $U$. Then $\phitrop((U')^{\rm an})=\phitrop(\Uan)$. 
\end{lem}

\proof 
Let $\omega \in \phitrop(\Uan)$. We note that $\varphi_{\rm trop}^{-1}(\omega)$ is a Laurent domain in  $\Uan$ and hence it has the same dimension as $U$. We conclude that $\varphi_{\rm trop}^{-1}(\omega)$ is not contained in the analytification of the lower dimensional Zariski-closed subset $U \setminus U'$ and hence $\omega \in \phitrop((U')^{\rm an})$. 
\qed

\begin{art} \rm \label{push-forward of cycles}
If $f: X_1 \rightarrow X_2$ is a morphism of varieties over $K$, then we define the {\it push-forward} of $X_1$ with respect to $f$ as the cycle $f_*(X_1):= \deg(f) \overline{f(X_1)}$, where the {\it degree} of $f$ is defined as $\deg(f):=[K(X_1):K(f(X_1))]$ if $f$ is generically finite and we set $\deg(f):=0$ if $[K(X_1):K(f(X_1))]=\infty$. By restriction, the push-forward can be defined in the same way on prime cycles of $X_1$ and extends by linearity to all cycles of $X_1$. 

Now let $\varphi:U \rightarrow T=\Tor$ be a moment map of the open subset $U$ of $X$. By \ref{tropical variety}, $$\Trop(\varphi_*(U)):=\deg(\varphi)\Trop(\overline{\varphi(U)})$$ is a tropical cycle  on $\rdop^r$. If $\varphi$ is generically finite, then this tropical cycle is of pure dimension $\dim(X)$ and the support is equal to $\phitrop(\Uan)$ (see Lemma \ref{density and tropicalization}). 
\end{art} 

The following result is called the {\it Sturmfels--Tevelev multiplicity formula}. It was proved by Sturmfels and Tevelev \cite{ST08} in the case of a trivial valuation and later generalized by Baker, Payne and Rabinoff \cite{BPR11} for every valued field.

\begin{prop} \label{Sturmfels-Tevelev multiplicity formula}
Let $\varphi':U' \rightarrow T'$ be a moment map of the non-empty open subset $U'$ of $X$ which refines the moment map $\varphi:U \rightarrow T$, i.e. there is an affine homomorphism $\psi:T' \rightarrow T$ such that $\varphi = \psi \circ \varphi'$ on $U' \subset U$. Then we have $$(\Trop(\psi))_*(\Trop(\varphi_*'(U')))=\Trop(\varphi_*(U))$$
in the sense of tropical cycles (see \ref{push-forward of tropical cycles}).
\end{prop}

\proof In fact, the Sturmfels--Tevelev multiplicity formula is the special case where $X=U'$ is a closed subvariety of $T'$ (see \cite{Gu12}, Theorem 13,17, for a proof in our setting deducing it from the original sources). In the general case, we conclude that 
$$(\Trop(\psi))_*(\Trop(\varphi_*'(U')))=\Trop(\psi_*((\varphi')_*(U')))=\Trop(\varphi_*(U')).$$
Since $U'$ is dense in $U$, the claim follows. \qed

\begin{art} \rm \label{final chart}
We will show that every open affine subset $U$ of $X$ has a canonical moment map. We note that the abelian group $M_U:=\Ocal(U)^\times/K^\times$ is free of finite rank (see \cite{Sa66}, Lemme 1). Here, we use that $K$ is algebraically closed (or at least that $X$ is geometrically reduced). We choose representatives  $\varphi_1, \dots , \varphi_r$ in $\Ocal(U)^\times$ of a basis. This leads to a moment map $\varphi_U:U \rightarrow T_U = \Spec(K[M_U])$. By construction, $\varphi_U$ refines every other moment map on $U$. Note that this moment map $\varphi_U$  is canonical up to (multiplicative) translation by an element of $T_U(K)$. 

Let $f:X' \rightarrow X$ be a morphism of algebraic varieties over $K$ and let $U'$ is an open subset of $X'$ with $f(U') \subset U$. Then $f^\sharp$ induces a homomorphism $M_U \rightarrow M_{U'}$ of lattices. We get a canonical affine homomorphism $\psi_{U,U'}:T_{U'} \rightarrow T_U$ of the canonical tori with $\psi_{U,U'} \circ \varphi_{U'} = \varphi_U \circ f$. This will be applied very often in the case where $U'$ is an open subset of $U$ in $X'=X$ and $f= \id$. Then we get a canonical affine homomorphism $\psi_{U,U'}:T_{U'} \rightarrow T_U$.
\end{art}

\begin{art} \rm \label{very affine}
Recall that an open subset $U$  of $X$ is called {\it very affine} if $U$ has a closed embedding into a multiplicative torus. Clearly, the following conditions are equivalent for an open affine subset $U$ of $X$: 
\begin{itemize}
 \item[(a)] $U$ is very affine;
 \item[(b)] $\Ocal(U)$ is generated as a $K$-algebra by $\Ocal(U)^\times$;
 \item[(c)] the canonical moment map $\varphi_U$ from \ref{final chart}   is a closed embedding.
\end{itemize}
The  intersection of  two very affine open subsets is again very affine (see the proof of Proposition \ref{chart lemma}). Moreover, the very affine open subsets of $X$ form a basis for the Zariski topology. We conclude that all local considerations can be done using very affine open subsets. 

On a very affine open subset, we will almost always use the canonical moment map $\varphi_U: U \rightarrow T_U$ which is a closed embedding by the above. To simplify the notation, we will set $\Trop(U)$ for the tropical variety of $U$ in $T_U$. It is a tropical cycle in $(N_U)_\rdop$, where $N_U$ is the dual abelian group of $M_U$. The tropicalization map will be denoted by $\trop_U:=(\varphi_U)_{\rm trop}:\Uan \rightarrow (N_U)_\rdop$. Recall that $\varphi_U$ is only determined up to translation by an element of $T_U(K)$ and hence $\trop_U$ and $\Trop(U)$ are only canonical up to an affine translation. This ambiguity is no problem as our constructions will be compatible with affine translations.
\end{art}

The following result of Ducros relates the local invariant $d(x)$ from \ref{local dimension} with tropical dimensions.

\begin{prop} \label{tropical stalk}
For $x \in \Xan$, there is a very affine  open neighbourhood $U$ of $x$ in $X$ such that for any open neighbourhood  $W$ of $x$ in the analytic topology of $\Uan$, there is a compact neighbourhood $V$ of $x$ in $W$ such that $\trop_U(V)$ is a finite union of $d(x)$-dimensional integral $\Gamma$-affine polytopes.
\end{prop}

\proof We choose rational functions $f_1, \dots, f_s$ on $X$ with $|f_1(x)|= \dots = |f_s(x)|=1$ such that the reductions $\widetilde{f_1}, \dots, \widetilde{f_s}$ form a transcendence basis of the residue field extension of $ {\mathscr H}(x)/K$. There are rational functions $g_1,\dots, g_t$ which are regular at $x$ such that $|g_1(x)|, \dots, |g_t(x)|$ form a basis of $(|{\mathscr H}(x)^\times|/|K^\times|) \otimes_\zdop \qdop$. By definition, we have $d(x)=s+t$. By (0.12) in \cite{Du12},  $f_1(x), \dots, f_s(x), g_1(x), \dots, g_t(x)$ reduce to a transcendence basis of the graded residue field extensions of ${\mathscr H}(x)/K$ in the sense of Temkin. There is a very affine open neighbourhood $U$ of $x$ in $X$ such that $f_1, \dots, f_s, g_1, \dots, g_t$ are invertible on $U$. Let $\varphi_1, \dots, \varphi_r \in \Ocal(U)^\times$ be the coordinates of the canonical moment map $\varphi_{U}:U \rightarrow T_{U}=\Tor$. Then the graded reductions of $\varphi_1, \dots, \varphi_r$ generate a graded subfield of the graded residue field extension of ${\mathscr H}(x)/K$. By construction, this graded subfield has transcendence degree $d(x)$ over the graded residue field of $K$. By  \cite{Du12}, Theorem 3.2, $\Trop_U(V)$ is a finite union of integral $\Gamma$-affine polytopes for every compact neighbourhood $V$ of $x$ in $\Uan$ which is strict in the sense of \cite{Be93}. For any open neighbourhood $W$ of $x$ in $\Uan$, Theorem 3.3 in \cite{Du12} shows that there is a compact strict neighbourhood $V$ of $x$ in $W$ such that $\trop_U(V)$ is a finite union of $d(x)$-dimensional  polytopes. \qed

\begin{art} \rm \label{tropical chart}
A {\it tropical chart} $(V,\varphi_U)$ on $\Xan$ consists of  an open subset $V$ of $\Xan$ contained in $\Uan$ for a very affine open subset $U$ of $X$ with $V=\trop_U^{-1}(\Omega)$ for some open subset $\Omega$ of  $\Trop(U)$. Here    the canonical moment map $\varphi_U:U \rightarrow T_U$ from \ref{final chart}  plays the role of (tropical) coordinates for $V$. By \ref{very affine}, $\varphi_U$ is an embedding. The condition $V=\trop_U^{-1}(\Omega)$ means that $V$ behaves well with respect to the tropical coordinates. In particular, $\trop_U(V)=\Omega$ is an open subset of $\Trop(U)$.

We say that the tropical chart $(V',\varphi_{U'})$ is a {\it tropical subchart} of $(V,\varphi_U)$ if $V'\subset V$ and $U' \subset U$. We note that the definition of tropical chart here is different from the tropical charts in \cite{CD12}, \S 3.1, which consist of an analytic morphism to a split torus and a finite union of polytopes containing the tropicalization. 
\end{art}

\begin{prop} \label{chart lemma}
The tropical charts on $\Xan$ have the following properties:
\begin{itemize}
 \item[(a)] They form a basis on $\Xan$, i.e. for every open subset $W$ of $\Xan$ and for every $x \in W$, there is a tropical chart $(V,\varphi_U)$ with $x \in V \subset W$. We may find such a $V$ such that the open subset $\trop_U(V)$  of $\Trop(U)$ is relatively compact.
 \item[(b)] The intersection $(V \cap V',\varphi_{U \cap U'})$ of tropical charts $(V,\varphi_U)$ and $(V',\varphi_{U'})$ is a tropical subchart of both.
 \item[(c)] If $(V,\varphi_U)$ is a tropical chart and if $U''$ is a very affine open subset of $U$ with $V \subset (U'')^{\rm an}$, then $(V,\varphi_{U''})$ is a tropical subchart of $(V,\varphi_U)$.
\end{itemize}
\end{prop}

\proof
To prove (a), we may assume that $X=\Spec(A)$ is a very affine scheme. A basis of $\Xan$ is formed by subsets of the form $V:=\{x \in X \mid s_1 <|f_1(x)|<r_1, \dots, s_k < |f_k(x)|<r_k\}$ with all $f_a \in A$ and real numbers $s_a < r_a$. Using the ultrametric triangle inequality as applied to $f_a+\pi$ for a non-zero $\pi \in K$ of small absolute value if $f_a(x)=0$, it is easy to see that we may choose the basis in such a way that $0<s_a$ for all $a=1, \dots k$. Note that $V$ is contained in the analytification of the very affine open subset $U:=\{x \in X \mid f_1(x) \neq 0, \dots , f_k(x) \neq 0\}$ of $X$.  It is obvious that $(V,\varphi_U)$ is a tropical chart proving (a).

To prove (b), let us consider the moment map
$$\Phi: U \cap U' \rightarrow T_U \times T_{U'}, \quad x \mapsto (\varphi_U(x), \varphi_{U'}(x)).$$
Since $X$ is separated, it is easy to see that $\Phi$ is a closed embedding and hence $U \cap U'$ is very affine.  
By definition of a tropical chart, $\Omega := \trop_U(V)$ (resp. $\Omega' := \trop_{U'}(V')$) is an open subset of $\Trop(U)$ (resp. $\Trop(U')$). Note that 
$$\Omega'':= \Phi_{\rm trop}((U \cap U')^{\rm an}) \cap (\Omega \times \Omega') \subset (N_U)_\rdop \times (N_{U'})_\rdop$$
is an open subset of $\Phi_{\rm trop}((U \cap U')^{\rm an})$. An easy diagram chase yields $\Phi_{\rm trop}^{-1}(\Omega'')=V \cap V'$. Since $\varphi_{U \cap U'}$ refines the moment map $\Phi$, we deduce that  $(V \cap V', \varphi_{U \cap U'})$ is a tropical chart. This proves (b).

Finally, we prove (c). Let $\psi:=\psi_{U,U''}:T_{U''} \rightarrow T_U$ be the canonical affine homomorphism from \ref{final chart}. Then we have $\trop_U=\Trop(\psi) \circ \trop_{U''}$ on $(U'')^{\rm an}$. 
Since $(V,\varphi_U)$ is a tropical chart, $\Omega:=\trop_U(V)$ is an open subset of $\Trop(U)$ and $V=\trop_U^{-1}(\Omega)$. 
Using $V \subset (U'')^{\rm an}$, we get $V= \trop_{U''}^{-1}(\Omega'')$ for the open subset $\Omega'':=\Trop(\psi)^{-1}(\Omega)$ of $\Trop(U'')$. We conclude that $(V,\varphi_{U''})$ is a tropical chart proving (c). \qed


\begin{rem} \rm \label{analytic generalizations}
In \cite{CD12}, everything is defined for an arbitrary analytic space. In Section \ref{analytic gen}, we will  compare their analytic constructions with our algebraic approach. 
\end{rem}

\section{Differential forms on algebraic varieties}

On a complex analytic manifold $M$, we use open analytic charts  $\varphi:U \rightarrow \cdop^r$ to define $(p,q)$-forms on $U$ by pull-back. The idea in the non-archimedean setting is similar replacing the above charts by tropical charts $(V,\varphi_U)$ from the previous section in order to pull-back Lagerberg's superforms to $\Uan$.

In this section, $K$ is an algebraically closed field endowed with a complete non-trivial non-archimedean absolute value $|\phantom{a}|$. Let $v:= - \log |\phantom{a}|$ be the associated valuation and let $\Gamma:=v(K^\times)$ be the value group. 
The theory could be done for arbitrary fields (see \cite{CD12}), but it is no serious restriction to assume that $K$ is algebraically closed as the theory is stable under base extension and in the classical setting, the analysis is also done over $\cdop$. We will introduce $(p,q)$-forms on the analytification $\Xan$ of a $n$-dimensional algebraic variety $X$ over $K$.

\begin{art} \rm \label{tropicalization}
We recall from \ref{tropical chart} that a tropical chart $(V,\varphi_U)$ consists of an open subset $V$ of $\Uan$ for a very affine open subset $U$ of $X$ such that $V=\trop_U^{-1}(\Omega)$ for an open subset $\Omega$ of $\Trop(U)$. Here,  $\varphi_U:U \rightarrow T_U$ is the canonical moment map. It is a closed embedding to the torus $T_U = \Spec(K[M_U])$. The tropical variety $\Trop(U)$ is a tropical cycle of $(N_U)_\rdop$ and $\trop_U:\Uan \rightarrow (N_U)_\rdop$ is the tropicalization map. The embedding $\varphi_U$ is only determined up to translation by an element in $T_U(K)$ and hence the tropical constructions are canonical up to integral $\Gamma$-affine isomorphisms.

Suppose that we have another tropical chart $(V',\varphi_{U'})$. Then $(V \cap V',\varphi_{U \cap U'})$ is a tropical chart (see Proposition \ref{chart lemma}) and we get a canonical affine homomorphism $\psi_{U, U \cap U'}:T_{U \cap U'} \rightarrow T_U$ of the underlying tori  with $\varphi_U= \psi_{U, U \cap U'}\circ \varphi_{U \cap U'}$ on $U \cap U'$ (see \ref{final chart}).  The associated  affine map $\Trop(\psi_{U, U \cap U'}):(N_{U \cap U'})_\rdop \rightarrow (N_U)_\rdop$ maps the tropical variety $\Trop(U \cap U')$ onto $\Trop(U)$ (use Lemma \ref{density and tropicalization}). 
Then we define the {\it restriction} of the superform $\alpha \in A^{p,q}(\trop_U(V))$ to a superform $\alpha|_{V \cap V'}$ on $\trop_{U \cap U'}(V \cap V')$ by using  the pull-back  to $\trop_{U \cap U'}(V \cap V')$  with respect to $\Trop(\psi_{U, U \cap U'})$. This plays a crucial role in the following definition:

\end{art}

\begin{Def} \rm \label{differential forms on varieties} 
A {\it differential form} $\alpha$ of bidegree $(p,q)$ on an open subset $V$ of $\Xan$ is given by a covering $(V_i)_{i \in I}$ of $V$ by tropical charts $(V_i,\varphi_{U_i})$ of $\Xan$ and superforms $\alpha_i \in A^{p,q}(\trop_{U_i}(V_i))$ such that $\alpha_i|_{V_i \cap V_j}=\alpha_j|_{V_i \cap V_j}$ for every $i,j \in I$.  If $\alpha'$ is another differential form of bidegree $(p,q)$ on $V$ given by $\alpha_j' \in A^{p,q}(\trop_{U'_j}(V_j'))$ with respect to the tropical charts $(V_j',\varphi_{U'_j})_{j \in J}$ covering $V$, then we consider $\alpha$ and $\alpha'$ as the same differential forms if and only if  $\alpha_i|_{V_i \cap V_j'}=\alpha_j'|_{V_i \cap V_j'}$ for every $i \in I$ and $j\in J$. We denote the space of $(p,q)$-differential forms on $V$ by $A^{p,q}(V)$. As usual, we define the space of differential forms on $V$ by $A(V):=\bigoplus_{p,q}A^{p,q}(V)$. The subspace of differential forms of degree $k \in \ndop$ is denoted by $A^k(V):=\bigoplus_{p+q =k}A^{p,q}(V)$.
\end{Def} 

\begin{art} \rm \label{properties of differential forms on varieties}
It is obvious from the definitions that the differential forms form a sheaf on $\Xan$. Using the corresponding constructions for superforms on tropical cycles, it is immediate to define the wedge product and differential operators $d$, $d'$, $d''$ on differential forms on $V$. By \ref{tropical variety}, we have $A^{p,q}(V)=\{0\}$ if $\max(p,q) > \dim(X)$.

For a morphism $\varphi:X' \rightarrow X$ and open subsets $V$ (resp. $V'$) of $\Xan$ (resp. $(X')^{\rm an}$) with $\varphi(V') \subset V$, we get a pull-back $\varphi^*:A^{p,q}(V) \rightarrow A^{p,q}(V')$ defined in the following way: Suppose that $\alpha \in A^{p,q}(V)$ is given by the covering $(V_i)_{i \in I}$ and the superforms  $\alpha_i  \in A^{p,q}(\trop_{U_i}(V_i))$  as above. Then there is a covering $(V_j')_{j \in J}$ of $V'$ by tropical charts $(V_j',\varphi_{U_j'})$ which is subordinate to $((\varphi^{\rm an})^{-1}(V_i))_{i \in I}$. This means that for every $j \in J$, there is $i(j) \in I$ with $V_j' \subset V_{i(j)}$ and   $\varphi(U_j') \subset U_{i(j)}$ for the corresponding very affine open subsets. Then $\varphi^*(\alpha)$ is the differential form on $V'$ given by the covering $(V_j')_{j \in J}$ and the superforms $\varphi^*(\alpha_{i(j)}) \in A^{p,q}(\Trop(U_j'))$. We leave the details to the reader. This construction is functorial as usual. 
\end{art}


\begin{rem} \rm \label{same superforms as in CD}
We obtain the same sheaf of differential forms on $\Xan$ as in \cite{CD12}, \S 3. In the latter reference, all analytic moment maps were used to define differential forms on $\Xan$ and so it is clear that our differential forms here are also differential forms in the sense of \cite{CD12}. To see the converse, we argue as follows: By Proposition \ref{chart lemma}, tropical charts $(V, \varphi_U)$ form a basis in $\Xan$. It follows from Proposition \ref{approximation} that an analytic moment map $\varphi:V \rightarrow (\Tor)^{\rm an}$ may be locally in $x \in V$ approximated by an algebraic moment map $\varphi':U' \rightarrow \Tor$ such that $(\varphi')_{\rm trop}=\trop \circ \varphi$ in an open neighbourhood of $x$ in $V$. Here, $U'$ is a suitable very affine open subset of $U$ with $x \in (U')^{\rm an}$. It follows from \cite{CD12}, Lemma 3.1.10, that we may use algebraic moment maps to define differential forms in the sense of \cite{CD12}. Using that $\varphi_{U'}$ factorizes through $\varphi'$ (see \ref{final chart}), we get the claim.
\end{rem}

\begin{Def} \rm \label{support of a differential form}
Let $\alpha$ be a differential form on an open subset $V$ of $\Xan$. The {\it support} of $\alpha$ is the complement in $V$ of the set of points $x$ of $V$ which have an open neighbourhood $V_x$ such that $\alpha|_{V_x}=0$. Let $A_c^{p,q}(V)$ be the space of differential forms of bidegree $(p,q)$ with compact support in $V$. 
\end{Def}

\begin{prop} \label{tropical forms}
Let $(V,\varphi_U)$ be a tropical chart of $\Xan$ and let $\alpha \in A^{p,q}(V)$ be given by $\alpha_U \in A^{p,q}(\trop_U(V))$. Then $\alpha =0$ in $A^{p,q}(V)$ if and only if $\alpha_U=0$ in $ A^{p,q}(\trop_U(V))$. 
\end{prop}

\proof See \cite{CD12}, Lemme 3.2.2. \qed

\begin{rem} \rm \label{tropical forms continued}
It follows from Proposition \ref{tropical forms}, that $\trop_U(\supp(\alpha))=\supp(\alpha_U)$ (see \cite{CD12}, Corollaire 3.2.3). Note however that not every differential form $\alpha$ on the tropical chart $(V, \varphi_U)$ is given by a single $\alpha_U \in A^{p,q}(\trop_U(V))$ as in Proposition \ref{tropical forms}.
\end{rem}

\begin{art} \rm \label{smooth functions}
In analogy with differential geometry on manifolds, we set $C^\infty(V):=A^{0,0}(V)$ for any open subset $V$ of $\Xan$ and a {\it smooth function} on $V$ is just a differential form of bidegree $(0,0)$. Since  tropicalization maps are continuous, it is clear that a smooth function is a continuous function on $V$. By the Stone-Weierstrass theorem, the space $C_c^\infty(V)$ of smooth functions with compact support in $V$ is a dense subalgebra of $C_c(V)$ (see \cite{CD12}, Proposition 3.3.5). 
\end{art}

\begin{Def} \rm \label{partition of unity}
Let $(V_i)_{i \in I}$ be an open covering of an open subset $V$ of $\Xan$. A   smooth {\it partition of unity} on $V$ with compact supports  subordinated to the covering $(V_i)_{i \in I}$ is a family $(\phi_j)_{j \in J}$ of non-negative smooth functions with compact support on $V$ with the following properties:
\begin{itemize}
 \item[(i)] The family $(\supp(\phi_j))_{j \in J}$ is locally finite on $V$.
 \item[(ii)] We have $\sum_{j \in J} \phi_j \equiv 1$ on $V$.
 \item[(iii)] For every $j \in J$, there is $i(j) \in I$ such that  $\supp(\phi_j) \subset V_{i(j)}$.
\end{itemize}
\end{Def}

\begin{prop} \label{existence of partition of unity}
Let $(V_i)_{ i \in I}$ be an open covering of an open subset $V$  of $\Xan$. Then there is a smooth partition of unity $(\phi_j)_{j \in J}$ on $V$ with compact supports subordinated to the covering $(V_i)_{i \in I}$. 
\end{prop}

\proof It is enough to show that for every $x \in V$, there is a non-negative smooth function $\phi$  with compact support in $V$ and with $\phi(x) >0$. Since $\Xan$ is a locally compact Hausdorff space which is also $\sigma$-compact, the open subset $V$ is paracompact and hence standard arguments from differential geometry yield the existence of the desired partition of unity (see \cite{Wa83}, Theorem 1.11).           



To prove the crucial claim at the beginning of the proof, we may assume that $V$ is coming from a tropical chart $(V,\varphi_U)$ 
(see Proposition \ref{chart lemma}). Then $\Omega := \trop_U(V)$ is a  open subset of $\Trop(U)$ with $\trop_U^{-1}(\Omega)=V$ and hence there is an open subset $\tilde{\Omega}$ in $(N_U)_\rdop$ with $\Omega = \tilde{\Omega} \cap \Trop(U)$. There is a smooth non-negative function $f$ on $(N_U)_\rdop$ with compact support in $\tilde{\Omega}$ such that $f(\trop_U(x))>0$. Since the tropicalization map is proper, the smooth function $\phi:=f \circ \trop_U$ has compact support in $V$ and hence $\phi$ fulfills the claim. \qed


\vspace{2mm} So far, we have seen properties of differential forms which are completely similar to the archimedean case. The next result of Chambert-Loir and Ducros (\cite{CD12}, Lemme 3.2.5) shows that the support of a differential form of degree at least one is disjoint from $X(K)$. 

\begin{lem} \label{support and d(x)}
Let $W$ be an open subset of $\Xan$. We consider $\alpha \in A^{p,q}(W)$ and  $x \in W$  with $d(x) < \max(p,q)$. Then $x \not \in \supp(\alpha)$.
\end{lem}

\proof Using Proposition \ref{chart lemma} and shrinking the open neighbourhood $W$ of $x$, we may assume that $W$ is a tropical chart $(W,\varphi_U)$ on which $\alpha$ is given by the superform $\alpha_U \in A^{p,q}(\trop_U(W))$. By Proposition \ref{tropical stalk}, there is a very affine open subset $U_x$ of $U$ and a compact  neighbourhood $V_x$ of $x$ in $(U_x)^{\rm an} \cap W$ such that $\trop_{U_x}(V_x)$ is of dimension $d(x)$. By Proposition \ref{chart lemma}, there is a tropical chart $(V',\varphi_{U'})$  with $x \in V' \subset V_x$ and $U' \subset U_x$. By \ref{final chart}, there is an affine homomorphism $\psi:T_{U'} \rightarrow T_U$ such that $\varphi_U = \varphi_{U'} \circ \psi$. Using the same factorization for the tropicalizations, we see that the restriction of $\alpha$ to $V'$ is given by $\Trop(\psi)^*(\alpha_U) \in  A^{p,q}(\trop_{U'}(V'))$. The inclusion $U_x \subset U$ yields that $\trop_{U_x}$ factorizes through $\trop_U$ (use \ref{final chart}). Since $V' \subset V_x$, we get  $\dim(\trop_U(V')) \leq \dim(\trop_{U_x}(V'))\leq d(x) <\max(p,q)$. As $\trop_U(V')=\Trop(\psi)( \trop_{U'}(V'))$, we conclude that $\Trop(\psi)^*(\alpha_U)=0$. This proves $\alpha=0$. \qed

\begin{cor} \label{support corollary}
Let $W$ be an open subset of $\Xan$ and let $U$ be a Zariski open subset of $X$. If $\alpha \in A^{p,q}(W)$ with $\dim(X \setminus U) < \max(p,q)$, then $\supp(\alpha) \subset W \cap \Uan$.
\end{cor}

\proof Let $x \in W \setminus \Uan$. Then \ref{local dimension} shows that $d(x) \leq \dim(X \setminus U) < \max(p,q)$. By Lemma  \ref{support and d(x)}, we get $x \not \in \supp(\alpha)$ proving the claim. \qed



\begin{prop} \label{single chart}
Let  $\alpha \in A^{p,q}_c(\Xan)$ be a differential form with $\max(p,q)=\dim(X)$. Then there is a very affine open subset $U$ of $X$ such that  $\supp(\alpha) \subset \Uan$ and such that $\alpha$ is given on $\Uan$ by a superform $\alpha_U \in A_c^{p,q}(\Trop(U))$. 
\end{prop}

\proof By assumption, the support of $\alpha$ is a compact subset of $\Xan$. We conclude that there are finitely many tropical charts $(V_i,\varphi_{U_i})_{i=1,\dots,s}$ covering $\supp(\alpha)$ such that $\alpha$ is given on $V_i$ by the superform $\alpha_i \in A^{p,q}(\trop_{U_i}(V_i))$. Recall that $\Omega_i:=\trop_{U_i}(V_i)$ is an open subset of $\Trop(U_i)$. By \ref{very affine}, $U:=U_1 \cap \dots \cap U_s$ is a non-empty very affine open subset of $X$. We define the open subset $V$ of $\Uan$ by $V:= \Uan \cap \bigcup_{i=1}^s V_i$. Since $\max(p,q)=\dim(X)$,  Corollary \ref{support corollary} yields $\supp(\alpha) \subset \Uan$.  Using \ref{final chart}, we see that $\trop_{U_i} =\Trop(\psi_i) \circ \trop_U$ for an affine homomorphism $\psi_i:T_U \rightarrow T_{U_i}$ of tori. Then we have  
$$\trop_U(V_i \cap \Uan)=(\Trop(\psi_i))^{-1}(\Omega_i) \cap \Trop(U)$$ and we denote this open subset of $\Trop(U)$ by $\Omega_i'$. It follows that the preimage of  $\Omega:= \bigcup_{i=1}^s \Omega_i'$  with respect to  $(\varphi_U)_{\rm trop}$ is equal to $V$. We conclude that $(V, \varphi_U)$ is a tropical chart of $\Xan$. Note that $\alpha$ is given on $\Uan \cap V_i$ by $\alpha_i':=\Trop(\psi_i)^*(\alpha_i) \in A^{p,q}(\Omega_i')$. By Proposition \ref{tropical forms}, $\alpha_i'$ agrees with $\alpha_j'$ on $\Omega_i' \cap \Omega_j'$ for every $i,j \in \{1, \dots, s\}$ and hence they define a superform $\alpha_U \in A^{p,q}(\Omega)$. By construction, $\alpha_U$ gives the differential form $\alpha$ on $V$. It follows from Remark \ref{tropical forms continued} that $\alpha_U$ has compact support in $\Omega$. Since $\alpha$ has compact support in $V$, we conclude that $\alpha_U$ is a superform on $\Trop(U)$ which defines $\alpha$ on $\Uan$. \qed

\begin{art} \rm \label{integration of differential forms on X}
Let  $\alpha \in A_c^{n,n}(W)$ for an open subset $W$ of $\Xan$, where  $n:=\dim(X)$. Obviously, we may view $\alpha$ as an $(n,n)$-form on $\Xan$ with compact support.  We call a very affine open subset $U$ as in Proposition \ref{single chart} a {\it very affine chart of integration for $\alpha$}. Then $\alpha$ is given by a superform $\alpha_U \in A_c^{n,n}(\Trop(U))$. We define the {\it integral of $\alpha$ over $W$} by
$$\int_W \alpha :=  \int_{\Trop(U)} \alpha_U.$$
Here, we view $\Trop(U)$ as a tropical cycle (see \ref{tropical variety}) and we integrate  as in \ref{tropical integration}. 
\end{art}

\begin{lem} \label{integration well defined}
For  $\alpha\in A_c^{n,n}(W)$, the following properties hold:
\begin{itemize}
\item[(a)] If $U$ is a very affine chart of integration for $\alpha$, then every non-empty very affine open subset $U'$ of $U$ is a very affine chart of integration for $\alpha$.
\item[(b)] The definition of $\int_W \alpha$ is independent of the choice of the very affine chart of integration for $\alpha$.
\end{itemize}
\end{lem}

\proof By Corollary \ref{support corollary}, $\supp(\alpha) \subset (U')^{\rm an}$ and (a) follows. To prove (b), it is enough to show
\begin{equation} \label{independence identity}
\int_{\Trop(U)} \alpha_U = \int_{\Trop(U')} \alpha_{U'}
\end{equation}
for a non-empty very affine open subset $U'$ of $U$ by using (a). 
  The differential form $\alpha$ is given on $\Uan$ (resp. $(U')^{\rm an}$) by $\alpha_U  \in A^{n,n}_c(\Trop(U))$ (resp. $ \alpha_{U'} \in A^{n,n}_c(\Trop(U'))$). By \ref{final chart}, there is an affine homomorphism $\psi:T_{U'} \rightarrow T_U$ of the underlying canonical tori such that $\varphi_U = \Trop(\psi) \circ \varphi_{U'}$. It follows that $\alpha$ is given on $U'$ also by $\Trop(\psi)^*(\alpha_U)$.  By Proposition \ref{tropical forms}, we have $\alpha_{U'}=  \Trop(\psi)^*(\alpha_U)$. The Sturmfels--Tevelev multiplicity formula  shows that $\Trop(\psi)_*(\Trop(U'))=\Trop(U)$ (see Proposition \ref{Sturmfels-Tevelev multiplicity formula}). Then Proposition  \ref{projection formula for tropical superforms} shows that
\eqref{independence identity} holds. \qed

\begin{prop} \label{linearity of integration}
Let $\lambda, \rho \in \rdop$ and let $\alpha, \beta \in A_c^{n,n}(W)$. Then we have
$$\int_W \lambda \alpha + \rho \beta = \lambda \int_W \alpha + \rho \int_W \beta.$$ 
\end{prop}

\proof By Lemma \ref{integration well defined}, we may choose a simultaneous very affine chart of integration for both $\alpha$ and $\beta$. Then the claim follows by the corresponding property of the integration of superforms. \qed

\vspace{2mm}

We have also {\it Stokes' theorem} for differential forms on the open subset $W$ of $\Xan$. Note that $W$ has trivial boundary in the algebraic situation \cite[Theorem 3.4.1]{Be90} and hence the boundary does not occur as in the version \cite[Theorem 3.12.1]{CD12} for analytic spaces.

\begin{thm} \label{Stokes for differential forms on X}
For $n:=\dim(X)$ and $\alpha \in A_c^{2n-1}(W)$, we have $\int_W d' \alpha = \int_W d'' \alpha =0$ and hence $\int_W d\alpha = 0$.
\end{thm}

\proof By Proposition \ref{single chart}, there is a  very affine open subset $U$ of $X$ such that  $\supp(\alpha) \subset \Uan$ and such that $\alpha$ is given on $\Uan$ by a superform $\alpha_U \in A^{2n-1}_c(\Trop(U))$. Then $U$ is a very affine chart of integration for $d'\alpha$ and $d''\alpha$ and the claim follows  from Proposition \ref{tropical Stokes} and Proposition \ref{tropical cycle iff closed}. \qed

\begin{rem} \rm \label{analytic integration}
Integration of differential forms on complex manifolds is defined by using a partition of unity with compact supports subordinated to a covering by holomorphic charts. Surprisingly, this was not necessary in our non-archimedean algebraic setting as we have  defined integration by using a single suitable tropical chart. In fact, the use of a smooth partition of unity $(\phi_j)_{j \in J}$ with compact supports subordinate to an open covering of $W$ by tropical charts $(V_i,\varphi_{U_i})_{i \in I}$ would not work here directly. To illustrate this, suppose that $\alpha \in A^{n,n}_c(W)$ is given on $V_i$ by $\alpha_i \in A^{n,n}(\trop_{U_i}(V_i))$. If the functions $\phi_j$ are of the form $\phi_j = f_j \circ \trop_{U_{i(j)}}$ for some $V_{i(j)} \supset \supp(\phi_j)$ and $f_j \in C^\infty_c(\trop_{U_{i(j)}}(V_{i(j)}))$, then we could set $\int_W \alpha = \sum_{j \in J} \int_{\Trop(U_{i(j)})} f_j \alpha_{i(j)}$. However, the functions $\phi_j$ could not be expected to have this form and so this approach fails.

Chambert-Loir and Ducros  define integration more generally for differential forms on paracompact good analytic spaces (see \cite{CD12}, \S 3.8). The idea is to use a covering by the interiors of affinoid subdomains.  Then there is a smooth partition of unity with  supports subordinated to this covering which reduces the problem to defining integration over an affinoid subdomain. But in the affinoid case, one can find a single tropical chart of integration similarly as in Proposition \ref{single chart}. It follows from Remark \ref{canonical calibration} and  Proposition \ref{comparison of tropical multiplicities} that both definitions give the same integral on the analytification of an algebraic variety.
\end{rem}

\section{Currents on algebraic varieties}

In this section, $K$ is an algebraically closed field endowed with a non-trivial non-archimedean complete absolute value $|\phantom{a}|$. We consider an open subset $W$ of $\Xan$ for an algebraic variety $X$ over $K$ of dimension $n$. Similarly as in the complex case, we will first define a topology on $A_c^{p,q}(W)$ and then we will 
 define currents as continuous linear functionals on this space. We will see that the Poincar\'e--Lelong equation holds for a rational function.

\begin{art} \rm \label{topology on the space of differential forms}
Let $(V_i, \varphi_{U_i})_{i \in I}$ be finitely many tropical charts contained in $W$ and let $\Delta_i$ be a polytope contained in  the open subset $\Omega_i:=\trop_{U_i}(V_i)$ of $\Trop(U_i)$. We consider the space $A^{p,q}(V_i,U_i,\Delta_i:i\in I)$ of $(p,q)$-forms $\alpha$ on $W$ with  support in $C:=\bigcup_{i\in I} \trop_{U_i}^{-1}(\Delta_i)$ such that $\alpha$ is given on $V_i$ by a superform $\alpha_i \in A^{p,q}(\Omega_i)$ for every $i \in I$. Since the tropicalization map is proper (see \ref{tropicalization map}), the set $C$ is compact. Similarly as in the complex case, we endow $A^{p,q}(V_i,U_i,\Delta_i:i\in I)$ with the structure of a locally convex space such that a sequence $\alpha_k$ converges to $\alpha$ if and only if all derivatives of the superforms $\alpha_{k,i}$ converge uniformly to the derivatives of the superform $\alpha_i$ on $\Delta_i$. Here, $\alpha_{k,i}$ (resp. $\alpha_i$) is the superform on $\Omega_i$ which defines $\alpha_k$ (resp. $\alpha$) on $V_i$ and we mean more precisely the derivatives of the coefficients of $\alpha_{k,i}|_{\Delta_i}$ (resp. $\alpha_i|_{\Delta_i}$). 
It follows easily from Proposition \ref{chart lemma} that $A_c^{p,q}(W)$ is the union of all spaces $A^{p,q}(V_i,U_i,\Delta_i:i\in I)$ with $(V_i,U_i,\Delta_i:i\in I)$ ranging over all possibilities as above. 
\end{art}

\begin{art} \rm \label{currents on varieties}
A {\it current} on an open subset $W$ of $\Xan$ is a linear functional $T$ on $A_c^{p,q}(W)$ such that the restriction of $T$ to all subspaces $A^{p,q}(V_i,U_i,\Delta_i:i\in I)$ is continuous. The space of currents is a $C^\infty(W)$-module denoted by $D_{p,q}(W)$. As usual (cf. \ref{supercurrents}), we define the differential operators  $d'$, $d''$ and $d:=d'+d''$ on the total space of currents $D(W):=\bigoplus_{p,q} D_{p,q}(W)$. Using partitions of unity from Proposition \ref{existence of partition of unity}, it is easy to show that the currents form a sheaf on $\Xan$. 
\end{art}

\begin{ex} \rm \label{measures}
A signed Radon measure $\mu$ on an open subset $W$ of $\Xan$ induces a current $[\mu] \in D_{0,0}(W)$ by setting $[\mu](f) := \int_\Xan f d\mu$ using ordinary integration theory on $\Xan$. Since the topology on $A_c^{0,0}(W)=C_c^\infty(W)$ is finer than the topology induced by the supremum norm, we conclude that $[\mu]$ is indeed a current on $W$. 
\end{ex}

\begin{rem} \rm \label{pushforward of currents}
Let $\varphi:X' \rightarrow X$ be a proper morphism of algebraic varieties over $K$. Then there is a linear map $\varphi_*:D_{p,q}((X')^{\rm an}) \rightarrow D_{p,q}(\Xan)$, where the {\it push-forward} $\varphi_*(T') \in D_{p,q}(\Xan)$ of $T' \in D_{p,q}((X')^{\rm an})$ is characterized by 
$$\varphi_*(T')(\alpha) = T'( \varphi^*(\alpha))$$
for every $\alpha \in A_c^{p,q}(\Xan)$. It follows from continuity of the map $\varphi^*:A_c^{p,q}(\Xan)) \rightarrow A_c^{p,q}((X')^{\rm an} )$ that $\varphi_*(T)$ is indeed a current on $T$. To define the push-forward, we need the fact that a proper algebraic morphism  induces a proper morphism between the analytifications which implies that the preimage of a compact subset in $\Xan$ is compact (see \cite{Be90}, Proposition 3.4.7). 
\end{rem}

\begin{ex} \rm \label{current of integration} 
We have the {\it current of integration} $\delta_X \in D_{2n}(\Xan)$ given by $\delta_X(\alpha) = \int_X \alpha$ for $\alpha \in A_c^{2n}(\Xan)$. More generally, we define the {\it current of integration along a closed $s$-dimensional subvariety $Y$ of $X$} as the push-forward of $\delta_Y \in D_{2s}(\Yan)$ to $\Xan$. By abuse of notation, we denote this element of $D_{2s}(\Xan)$ also by $\delta_Y$. By linearity in the components, we define the current of integration along a cycle on $X$. If $W$ is an open subset of $\Xan$, then we get a current $\delta_W \in D_{2n}(W)$ by restricting $\delta_X$. 
\end{ex}

\begin{art} \rm \label{wedge product of currents with differential forms}
Let $T \in D_{p,q}(W)$ and  $\omega \in A^{r,s}(W)$ for an open subset $W$ of $\Xan$. Then we define $T \wedge \omega \in D_{p-r, q-s}(W)$ by $(T \wedge \omega)(\alpha) = T(\omega \wedge \alpha)$ for $\alpha \in A_c^{p-r,q-s}(W)$. Since the wedge product with a given form is a continuous operation on $A_c(W)$, it is clear that $T \wedge \omega$ is really a current on $W$.
\end{art}

\begin{ex} \rm \label{current associated to differential forms on varieties}
For $\omega \in A^{r,s}(W)$, the {\it current $[\omega] \in D_{n-r,n-s}(W)$ associated to $\omega$} is defined by $[\omega]:=\delta_W \wedge \omega$ and we get an injective linear map $a:A^{r,s}(W) \rightarrow D_{n-r,n-s}(W)$ given by $a(\omega):=[\omega]$.
\end{ex}

\begin{prop} \label{measures for topdimensional forms}
Let $\omega \in A^{2n}(W)$ for an open subset $W$ of $\Xan$. Then there is a unique signed Radon measure $\mu$ on $W$ such that $\int_W f d\mu = [\omega](f)$ for every $f \in C_c^\infty(W)$. If $\omega$ has compact support, then we have $|\mu|(W) < \infty$. 
\end{prop}

\proof It is easy to prove that $[\omega]$ induces a continuous linear functional on  $C_c^\infty(W)$ where this locally convex vector space is endowed with the subspace topology of $C_c(W)$. By  \ref{smooth functions}, this subspace is dense and hence the Riesz representation theorem proves the first claim. If $\omega$ has compact support, then $\supp(\mu)$ is also compact and the last claim follows. \qed

\begin{art} \rm \label{locally integrable differential forms on varieties}
Let us again consider an open subset $W$ of $\Xan$. A  function $f:W \rightarrow \rdop \cup \{\pm \infty\}$ is called {\it locally integrable} if $f$ is integrable with respect to the measure $\mu$ associated to any $\omega \in A_c^{2n}(W)$. Then we write $\int_W f \omega := \int_W f d\mu$. 

For a locally integrable function $f$ on $W$ and $\eta \in A^{p,q}(W)$, we define $[f \cdot \eta] \in D_{p,q}(W)$ by 
$[f \cdot \eta](\alpha):=\int_W f \eta \wedge \alpha$ for every $\alpha \in A_c^{n-p,n-q}(W)$.
\end{art}

Chambert-Loir and Ducros proved the {\it Poincar\'e--Lelong equation} for rational functions:

\begin{prop} \label{Poincare-Lelong for rational functions}
Let $f$ be a rational function on $X$ which is not identically zero. Then $\log|f|$ is a locally integrable function on $\Xan$ and we have $d'd'' [\log|f|] = \delta_{\Div(f)}$. 
\end{prop}

\proof See \cite{CD12}, Theorem 4.6.5. \qed

\section{Generalizations to analytic spaces} \label{analytic gen}

The final section shows how our notions fit with the paper \cite{CD12}. While we restricted to the algebraic case, the paper of Chambert-Loir and Ducros works for arbitrary analytic spaces. We assume that the reader is familiar with the theory of analytic spaces as given in \cite{Be93} or \cite{Te10}. For simplicity, we assume again that $K$ is algebraically closed, endowed with a non-trivial non-archimedean complete absolute value $|\phantom{a}|$ with corresponding valuation $v:=-\log|\phantom{a}|$ and that all occurring analytic spaces are strict in the sense of \cite{Be93}. This situation can always be obtained by base change without changing the theory of differential forms and currents. As usual, we use the value group $\Gamma:=v(K^\times)$.

\begin{art} \rm \label{tropicalization moment map}
Let $Z$ be a compact analytic space over $K$. An {\it analytic moment map} on $Z$ is an analytic morphism $\varphi:Z \rightarrow \Tan$ for a split torus $T=\Tor$ over $K$ as before. Let $M$ be the character group of $T$, then we have $T=\Spec(K[M])$. The map $\phitrop := \trop \circ \varphi:Z \rightarrow N_\rdop$ is called the {\it tropicalization map} of $\varphi$ and we may use the coordinates on $T$ to identify $N_\rdop$ with $\rdop^r$.
\end{art}

 The next result shows that for the construction of differential forms in the algebraic case, we may restrict our attention to algebraic moment maps.

\begin{prop} \label{approximation}
Let $X$ be an algebraic variety over $K$ and let $\varphi:W \rightarrow \Tan$ be an analytic moment map defined on an open subset $W$ of $\Xan$. For every $x \in W$, there is a very affine open subset $U$ of $X$  with an algebraic moment map  $\varphi':U \rightarrow T$ and an open neighbourhood $V$ of $x$ in $\Uan \cap W$ such that $\phitrop=\varphi'_{\rm trop}$ on $V$. 
\end{prop}

\proof We may assume that $X = \Spec(A)$. Similary as in the proof of Proposition \ref{chart lemma}, there is a  neighbourhood $V':=\{x \in X \mid s_1 \leq |f_1(x)| \leq r_1, \dots, s_k \leq |f_k(x)| \leq r_k\}$ of $x$ in $W$ with all $f_a \in A$ and real numbers $0 <s_a < r_a$. We may assume that $f_1, \dots, f_k$ form an affine coordinate system $y_1, \dots, y_k$ on $X$.  Using coordinates on $T=\Tor$, the moment map $\varphi$ is given by analytic functions $\varphi_1, \dots, \varphi_r$ on $W$ which restrict to strictly convergent Laurent series in $y_1, \dots, y_k$ on $V'$. Cutting the Laurent series in sufficiently high positive and negative degree, we get Laurent polynomials  $p_1, \dots, p_r$ with $|p_a| = |\varphi_a|$ on $V'$ for $a=1, \dots, r$. By Proposition \ref{chart lemma}, there is a very affine open subset $U$ of $X$ such that $\Uan$ contains $x$ and such that $p_1, \dots p_r$ define an algebraic moment map $\varphi':U \rightarrow  T$ with $\phitrop=\varphi'_{\rm trop}$ on $V'$. Choosing a neighbourhood $V$ of $x$ in $\Uan \cap V'$, we get the claim. \qed


\vspace{2mm}

We have the following generalization of the Bieri--Groves theorem. 
Working with analytic spaces, the  boundary $\partial Z$ of $Z$ becomes an issue.

\begin{thm}[Berkovich, Ducros] \label{analytic Bieri-Groves}
If $Z$ is a compact analytic space over $K$ of dimension $n$ and if $\varphi:Z \rightarrow \Tan$ is an analytic moment map, then $\phitrop(Z)$ is a finite union of integral $\Gamma$-affine polytopes of dimension at most $n$. Moreover, $\phitrop(\partial Z)$ is contained in a finite union of integral $\Gamma$-affine polytopes of dimension $\leq n-1$. If $Z$ is affinoid, then $\phitrop(\partial Z)$ is equal to a finite union of such polytopes.
\end{thm}

\proof The first claim is due to Berkovich and the remaining claims are due to Ducros  (see \cite{Du12}, Theorem 3.2). \qed

\begin{art} \label{analytic tropical variety} \rm 
We consider now a compact analytic space $Z$  over $K$ of pure dimension $n$. Theorem \ref{analytic Bieri-Groves} shows that  the {\it tropical variety} $\phitrop(Z)$ is the support of an integral $\Gamma$-affine polytopal complex in $N_\rdop$. Our next goal is to endow this complex with canonical tropical multiplicities. This will lead to the definition of a weighted polytopal complex $(\phitrop)_*(\cyc(Z))$ which is canonical up to subdivision. 

If $\dim(\phitrop(Z))<n$, then we set $(\phitrop)_*(\cyc(Z))=0$ meaning that we choose all tropical weights equal to zero. It remains to consider the case $\dim(\phitrop(Z))=n$. 
We choose a generic surjective homomorphism $q:T\rightarrow T'$ onto a split multiplicative torus $T'=\Spec(K[M'])$ of rank $n=\dim(Z)$. Generic means that the corresponding linear map $F:=\Trop(q)$ is injective on every polytope contained in $\phitrop(Z)$. By Theorem \ref{analytic Bieri-Groves}, there is an integral $\Gamma$-affine polytopal complex $\Ccal$ in $N_\rdop$ with $|\Ccal|=\phitrop(Z)$ such that $F(\tau)$ is disjoint from $(q \circ \varphi)_{\rm trop}(\partial Z)$ for every $n$-dimensional  face $\sigma$ of $\Ccal$ and $\tau := \relint(\sigma)$. By passing to a subdivision, we may assume that $F_*(\Ccal)$ is a polyhedral complex in $N_\rdop'$ as in \ref{push-forward of tropical cycles}, where $N'$ is the dual of $M'$ as usual.

We identify $F(\tau) \subset N'_\rdop \cong \rdop^{n}$ with an open subset of the skeleton $S((T')^{\rm an})$ as in Remark \ref{tropicalization as a retraction}. Then it is clear that $q \circ \varphi$ restricts to a map $(q \circ \varphi)^{-1}(\tau) \rightarrow F(\tau)$ which agrees with $F \circ \phitrop$  using the identification $S((T')^{\rm an})=N'_\rdop$. It is shown in \cite{CD12}, \S 2.4, that this restriction of $q \circ \varphi$  is a finite flat and surjective morphism which means that   every point $p$ of $F(\tau)$ has a neighbourhood $W'$ in $(T')^{\rm an}$ such that $(q \circ \varphi)^{-1}(W') \rightarrow W'$ has these properties.  Using that $F^{-1}(F(\tau))=\coprod_{\tau'} \tau'$,  where  $\tau'$ is ranging over all open faces of $\Ccal$ with $F(\tau')=F(\tau)$, we get 
$$(q \circ \varphi)^{-1}(F(\tau))= \coprod_{\tau'} \varphi_{\rm trop}^{-1}(\tau') \cap (q \circ \varphi)^{-1}(F(\tau)).$$
We conclude that the map $\varphi_{\rm trop}^{-1}(\tau) \cap (q \circ \varphi)^{-1}(F(\tau))\rightarrow F(\tau)$ is finite, flat and surjective. Again, this has to be understood in some open neighbourhoods. Since $\tau$ is connected, the  corresponding degree depends only on $\tau$ and not on the choice of $p$. We denote this degree by $[\varphi_{\rm trop}^{-1}(\tau) : F(\tau)]$.

Recall that $N_\sigma$ is the canonical lattice in the affine space generated by $\sigma$. Then the character lattice $M'$ of $T'$ is of finite index in $M_\sigma = \Hom(N_\sigma,\zdop)$.
\end{art}

\begin{Def} \rm \label{analytic tropical multiplicity}
Using the notation from above, the {\it tropical multiplicity} $m_\sigma$ along $\sigma$ is defined by 
$$m_\sigma:=[\varphi_{\rm trop}^{-1}(\tau) : F(\tau)]\cdot [M_\sigma : M']^{-1}.$$
Furthermore, $(\phitrop)_*(\cyc(Z))$ is the weighted polyhedral complex $\Ccal$ endowed with these tropical multiplicities. The weights might be rational numbers, at least we have no argument that they are integers in the analytic case. \end{Def}

\begin{rem} \rm \label{canonical calibration}
It is not so easy to show that the tropical multiplicity is well-defined, i.e. independent of the choice of $q$. Chambert-Loir and Ducros do not use tropical multiplicities, but the latter are equivalent to the canonical calibration introduced in \cite{CD12}, \S 3.5. To summarize this construction, let $e_1, \dots, e_n$ (resp. $f_1, \dots, f_n$) be a basis of $M'$ (resp. $M_\sigma$). Then the canonical calibration of $\sigma$ is defined as 
$$[\varphi_{\rm trop}^{-1}(\tau) : F(\tau)] \cdot (F|_{(N_\sigma)_\rdop})^*(e_1 \wedge \dots \wedge e_n) \in \Lambda^n((N_\sigma)_\rdop)$$
together with the orientation induced by the pull-back of $e_1, \dots, e_n$ with respect to the linear isomorphism $F|_{(N_\sigma)_\rdop}$. The canonical calibration is equal to the calibration $m_\sigma f_1 \wedge \dots \wedge f_n$ together with the orientation induced by $f_1, \dots , f_n$. Since the canonical calibration does not depend on the choice of $q$ up to refinement (\cite{CD12}, \S 3.5), the same is true for the tropical multiplicities.
\end{rem}

\begin{rem} \rm \label{cycle of analytic space}
One can define the irreducible components of an analytic space (see \cite{Con99}). A compact analytic space $Z$ has finitely many irreducible components $Z_i$. Then we  define the cycle $\cyc(Z)$ associated to $Z$ as a positive formal $\zdop$-linear combination of the irreducible components $Z_i$ by restriction to affinoid subdomains and then by glueing (see \cite{Gu98}, \S 2). One can show that the weighted $n$-dimensional polyhedral complex $(\phitrop)_*(\cyc(Z))$ depends only on $\cyc(Z)$ and this dependence is linear. We leave the details to the reader.
\end{rem}

The next result shows that the Sturmfels--Tevelev multiplicity formula holds for analytic spaces.

\begin{prop} \label{analytic Sturmfels-Tevelev}
Let $Z$ be a compact analytic space over $K$ of pure dimension $n$, let $\varphi:Z \rightarrow \Tan$ be an analytic moment map and let $\psi:T \rightarrow T'$ be an affine homomorphism of tori. Then we have $$\Trop(\psi)_*((\phitrop)_*(\cyc(Z)) )= ((\psi \circ \varphi)_{\rm trop})_*(\cyc(Z)).$$
 \end{prop}

\proof The corresponding statement for canonical calibrations is shown in \cite{CD12}, Lemma 3.5.2, and hence the claim follows from Remark \ref{canonical calibration}. \qed

\begin{prop} \label{analytic balancing}
Let $Z$ be a compact analytic space over $K$ of pure dimension $n$ and let $\Ccal$ be the same integral $\Gamma$-affine polytopal complex  with support  $\phitrop(Z)$ as in \ref{analytic tropical variety}. Then for every $(n-1)$-dimensional polyhedron $\rho$ of $\Ccal$ not contained in $\phitrop(\partial Z)$, the balancing condition 
$$\sum_{\sigma \in \Ccal_n, \,\sigma \supset \rho} m_\sigma \omega_{\rho, \sigma} \in N_\rho$$
from \ref{tropical cycle} holds in $\rho$. 
\end{prop}

\proof Chambert-Loir and Ducros prove in \cite{CD12}, Theorem 3.6.1, that $\rho$ is harmonious in $\Ccal$ which is a condition for the canonical calibration  equivalent to the balancing condition  by Remark \ref{canonical calibration}. \qed

\begin{art} \label{common setup} \rm
In an algebraic setting, our goal is to compare the tropical multiplicities introduced in \ref{tropical weights} with the ones from Definition \ref{analytic tropical multiplicity}. Let us consider an algebraic variety $X$ over $K$ of dimension $n$ and an algebraic moment map $\varphi:X\rightarrow T=\Spec(K[M])\cong \Tor$ over $K$.  
Note that $\phitrop(\Xan)=\Trop(\overline{\varphi(X)})$. 
We endow $\phitrop(\Xan)$ 
with the tropical multiplicities $m_\sigma^{\rm alg}$ of the tropical cycle $\Trop(\varphi_*(X)):=\deg(\varphi)\Trop(\overline{\varphi(X)})$ of $N_\rdop$. 

The analytification $\Xan$ is not compact (unless  $n=0$), but as $\partial X = \emptyset$, we can define tropical multiplicities in the same analytic manner as in Definition \ref{analytic tropical multiplicity}. This means that we choose a generic projection $q:T \rightarrow T'=\Spec(K[M'])$ onto a torus $T'$ of rank $n$ and an integral $\Gamma$-affine polyhedral complex $\Ccal$ with support equal to $\phitrop(\Xan)$ such that $F_*(\Ccal)$ is a polyhedral complex on $N'_\rdop$ for the associated linear map $F:N_\rdop \rightarrow N'_\rdop$. For every $\sigma \in \Ccal_n$ and $\tau := \relint(\sigma)$, we define $$m_\sigma^{\rm an}:=[\varphi_{\rm trop}^{-1}(\tau) : F(\tau)]\cdot [M_\sigma : M']^{-1}$$
as in Definition \ref{analytic tropical multiplicity}. Since the tropical multiplicities $m^{\rm alg}$ and $m^{\rm an}$ are  compatible with subdivision, we may assume that the underlying integral $\Gamma$-affine polyhedral complex $\Ccal$ is the same in both definitions. 
\end{art}

Now we are ready to compare these two tropical multiplicities. 

\begin{prop} \label{comparison of tropical multiplicities}
Let $\varphi:X \rightarrow T$ be an algebraic moment map and $n:=\dim(X)$. Using the notations from above, we have $m_\sigma^{\rm an}= m_\sigma^{\rm alg}$ for every $\sigma \in \Ccal_n$. 
\end{prop}

\proof The following argument is quite close to the proof of the Sturmfels--Tevelev formula given by Baker, Payne and Rabinoff (see \cite{BPR11}, Theorem 8.2). We may assume that $\varphi$ is generically finite, otherwise $\phitrop(\Xan)$ has dimension $<n$ and all tropical multiplicities are zero. Let $Y$ be the closure of $\varphi(X)$ in $T$ and let $q:T \rightarrow T'$ be a generic homomorphism onto a split torus $T'=\Spec(K[M'])$ of rank $n$ with associated linear map $F:N_\rdop \rightarrow N'_\rdop$. There is an open dense subset $U$ of $Y$ such that $\varphi$ is finite over $U$. Since the tropical multiplicities $m^{\rm an}$ and $m^{\rm alg}$ are compatible with subdivision of the polyhedral complex $\Ccal$, we may assume that $\Trop(Y \setminus U)$ is contained in the support of $\Ccal_{n-1}$. 


Let $Y'$ be the closure of $q(Y)$ in $T'$ and let $\omega \in \tau \cap N_\Gamma$. 
We consider the affinoid subdomains $U_\omega := \trop^{-1}(\omega)$ in $\Tan$ and $U'_{\omega'}:=\trop^{-1}(\omega')$ in $(T')^{\rm an}$. By finiteness of $\varphi$ over $U$, the set $X_\omega:=(\varphi^{\rm an})^{-1}(U_\omega)=\varphi_{\rm trop}^{-1}(\omega)$ is an affinoid subdomain of $\Xan$ and $\varphi$ restricts to a finite morphism $X_\omega \rightarrow Y_\omega := \Yan \cap U_\omega$. Let $\Xcal_\omega, \Ycal_\omega, \Ucal_\omega, \Ucal'_{\omega'}$ be the canonical formal affine $\kcirc$-models of $X_\omega, Y_\omega, U_\omega, U'_{\omega'}$ associated to the algebra of power bounded elements in the corresponding affinoid algebra. Moreover, let $\overline{Y_\omega}$ be the closure of $Y_\omega$ in $\Ucal_\omega$. Then we have canonical  morphisms
\begin{equation} \label{three morphisms}
 \Xcal_\omega \stackrel{\varphi}{\rightarrow} \Ycal_\omega \stackrel{\iota}\rightarrow \overline{Y_\omega} \stackrel{q}{\rightarrow} \Ucal'_{\omega'}
\end{equation}
of admissible formal affine schemes over $\kcirc$ in the sense of Bosch, L\"utkebohmert and Raynaud (see \cite{BL93}, \S 1). We claim that all these morphisms are finite and surjective. Obviously, the generic fibres of the first and second morphism are finite and surjective. To see that the generic fibre of the third morphism is finite, we note first that $F^{-1}(\omega') \cap \Trop(Y)$ is finite by construction of $q$ and hence $q^{-1}(U'_{\omega'}) \cap \Yan$ is in the relative interior of an affinoid subdomain of $\Tan$ which is contained in $q^{-1}(U_{\omega'}')$. We conclude that $q^{-1}(U'_{\omega'}) \cap \Yan \rightarrow U'_{\omega'}$ is a proper map (see the proof of Theorem 4.31 in \cite{BPR11} for more details about the argument). Since $q^{-1}(U'_{\omega'}) \cap \Yan$ is the  disjoint union of the finitely many affinoids $U_\rho \cap \Yan$, $\rho \in F^{-1}(\omega') \cap \Trop(Y)$, we conclude that $q$ induces a proper morphism $Y_\omega \rightarrow U'_{\omega'}$ of affinoids. By Kiehl's direct image theorem (\cite{BGR84}, Theorem 9.6.3/1), this morphism is finite and hence also surjective using dimensionality arguments. We conclude that all three morphisms in \eqref{three morphisms} are surjective and finiteness follows from  \cite{BPR11}, Proposition 3.13. 

The degree  $[X_\omega : U'_{\omega'}]$ of $X_\omega$ over the affinoid torus $U'_{\omega'}$ is well-defined as $U'_{\omega'}$ is irreducible (see \cite{BPR11}, Section 3, for a discussion of degrees). Since the degree does not change by passing to an affinoid subdomain of $U'_{\omega'}$ (see \cite{BPR11}, Proposition 3.30), we get
\begin{equation} \label{a}
[\varphi_{\rm trop}^{-1}(\tau) : F(\tau)]= [X_\omega : U'_{\omega'}].
\end{equation}
The projection formula (\cite{BPR11}, Proposition 3.32) shows
\begin{equation} \label{b}
[X_\omega : U'_{\omega'}]= \sum_B [B: (\Ucal'_{\omega'})_s] = \sum_B [B: ({\mathbb G}_m^{n})_{\ktilde}],
\end{equation}
 where $B$ ranges over all irreducible components of $(\Xcal_\omega)_s$. We conclude from \eqref{a} and \eqref{b} that
\begin{equation} \label{c}
[\varphi_{\rm trop}^{-1}(\tau) : F(\tau)]= \sum_C \sum_{\text{$B$ over $C$}} [B:C] \cdot [C: ({\mathbb G}_m^{n})_{\ktilde}],
\end{equation}
where $C$ ranges over all irreducible components of $(\overline{Y_\omega})_s$ and $B$ ranges over all irreducible components of $(\Xcal_\omega)_s$ mapping onto $C$. Since the special fibre of $\overline{Y_\omega}$ is isomorphic to the initial degeneration $\In_\omega(Y)$, all irreducible components $C$ are isomorphic to the torus $\Spec(\ktilde[M_\sigma])$ (see \cite{BPR11}, Theorem 4.29) proving
\begin{equation} \label{d}
[C: ({\mathbb G}_m^{n})_{\ktilde}]= [M_\sigma: M'].
\end{equation}
Using \eqref{c} and \eqref{d}, we get
\begin{equation} \label{e}
 m_\sigma^{\rm an}= [\varphi_{\rm trop}^{-1}(\tau) : F(\tau)] \cdot [M_\sigma : M']^{-1} = \sum_C \sum_{\text{$B$ over $C$}} [B:C].
\end{equation}
Since $X_\omega$ is the preimage of the affinoid subdomain $Y_\omega$ of $\Tan$, we deduce from \cite{BPR11}, Proposition 3.30, that $X_\omega$ is of  degree $\deg(\varphi)$ over $Y_\omega$  and hence the projection formula again shows the equality
\begin{equation} \label{f}
\deg(\varphi) \cyc((\overline{Y_\omega})_s)= (\iota \circ \varphi)_*( \cyc((\Xcal_\omega)_s)
\end{equation}
of cycles in $(\Ucal_\omega)_s$. Inserting \eqref{f} in \eqref{e} by using that the special fibre of $\Xcal_\omega$ is reduced, we get
$$m_\sigma^{\rm an}=\deg(\varphi) \sum_C m(C,(\overline{Y_\omega})_s),$$
where $m(C,(\overline{Y_\omega})_s)$ is the multiplicity of the irreducible component $C$ in the special fibre of $\overline{Y_\omega}$. By definition, the right hand side is equal to $m_\sigma^{\rm alg}$ which proves the claim.
\qed

\begin{rem} \rm \label{new proof of balancing}
Note that in the algebraic case, Proposition \ref{comparison of tropical multiplicities} yields that the tropical multiplicities in Definition \ref{analytic tropical multiplicity} are well-defined integers, i.e. independent of the choice of the generic projection $q$. Moreover, the argument of Chambert-Loir and Ducros for Proposition \ref{analytic balancing} gives a new proof for the classical balancing condition for tropical varieties which is based mainly on degree considerations.
\end{rem}

{\small Walter Gubler, Universit\"at Regensburg, Fakult\"at f\"ur Mathematik,   Universit\"atsstrasse 31, 
 D-93040 Regensburg, walter.gubler@mathematik.uni-regensburg.de}

\end{document}